\documentclass[preprint, 11pt]{elsarticle}

\usepackage[T1]{fontenc}
\usepackage[utf8]{inputenc}
\usepackage[hmargin=0.9in]{geometry}
\usepackage[babel,german=quotes]{csquotes}
\usepackage{subcaption}
\usepackage{graphicx}
\usepackage[colorlinks=true,citecolor=blue,urlcolor=blue]{hyperref}
\usepackage{caption}
\usepackage{amsmath}
\usepackage{amsthm}
\usepackage{amssymb}
\usepackage{bm, amsfonts}
\usepackage{xcolor,soul}
\usepackage{fixmath}
\usepackage{tikz}
\usepackage{bm}
\usepackage{titlesec}
\usepackage{multirow}
\usepackage{indentfirst}
\usepackage{mathtools}
\usepackage{float}
\usepackage[normalem]{ulem}
\usepackage{comment}
\useunder{\uline}{\ul}{}

\newcommand{\diff}{\mathrm{d} }
\newcommand{\R}{\mathbb{R}}
\newcommand{\N}{\mathcal{N}}

\newcommand{\x}{\mathbf{x}}
\newcommand{\s}{\mathbf{s}}
\newcommand{\vel}{\mathbf{v}}
\newcommand{\n}{\mathbf{n}}
\newcommand{\f}{\mathbf{f}}

\newcommand{\con}{\mathbf{c}}

\newcommand{\D}{\mathcal{D}}
\newcommand{\dv}{\nabla\cdot}
\newcommand{\Sp}{\mathbb{S}^{d-1}}
\newcommand{\OO}{\bm{\Omega}}
\newcommand{\pz}{\psi^{(0)}}
\newcommand{\po}{\boldsymbol{\psi}^{(1)}}
\newcommand{\pt}{\boldsymbol{\psi}^{(2)}}

\newcommand{\dx}{\,\diff\x}

\newcommand{\mev}{\,\mathrm{MeV}}
\newcommand{\cm}{\,\mathrm{cm}}

\usepackage{xcolor}

\colorlet{shadecolor}{green}

\titleformat{\paragraph}
{\normalfont\normalsize\itshape}{\theparagraph}{1em}{}
\titlespacing*{\paragraph}
{0pt}{3.25ex plus 1ex minus .2ex}{1.5ex plus .2ex}

\newdefinition{remark}{Remark}

\makeatletter
\def\ps@pprintTitle{%
\let\@oddhead\@empty
\let\@evenhead\@empty
\def\@oddfoot{
\footnotesize\itshape
\ifx\@journal\@empty Elsevier
\else\@journal\fi
\hfill\today
}%
\let\@evenfoot\@oddfoot}
\makeatother

\begin{document}

\begin{frontmatter}
\title{Realizability-preserving finite element discretizations of the $M_1$ model for dose calculation in proton therapy}
\author[TUD]{Paul Moujaes\corref{cor1}}
\ead{paul.moujaes@math.tu-dortmund.de}

\author[TUD]{Dmitri Kuzmin}
\ead{kuzmin@math.uni-dortmund.de}

\address[TUD]{Institute of Applied Mathematics (LS III), TU Dortmund University\\ Vogelpothsweg 87,
  44227 Dortmund, Germany}

\author[WPE,WTZ,DKTK,TUDOPhy]{Christian Bäumer}

\address[WPE]{West German Proton Therapy Centre Essen
(WPE) gGmbH\\ Am Mühlenbach 1, 45147 Essen, Germany}
\address[WTZ]{West German Cancer Center (WTZ), Hufelandstr. 55, 45147 Essen, University Hospital Essen, Essen, Germany}
\address[DKTK]{German Cancer Consortium (DKTK), Hufelandstr. 55, 45147 Essen, Germany }
\address[TUDOPhy]{Department of Physics, TU Dortmund University, Otto-Hahn-Str. 4, 44227 Dortmund, Germany}
\ead{Christian.Baeumer@uk-essen.de}

\cortext[cor1]{Corresponding author}

\journal{}

\begin{abstract}
  We present a deterministic framework for proton therapy dose calculation based on finite element discretizations of the energy-dependent $M_1$ moment model. The nonlinear $M_1$ system is derived from the Fokker--Planck equation for charged particles and closed using an entropy-based approximation of the second moment. Energy is treated as a pseudo-time coordinate. The zeroth and first moments of the proton fluence are evolved backward in energy. To ensure hyperbolicity and physical admissibility, we employ a monolithic convex limiting (MCL) strategy. Representing the standard continuous Galerkin discretization in terms of auxiliary `bar' states, we construct a nonlinear scheme that is provably invariant domain preserving (IDP) w.r.t. convex realizable sets consisting of all admissible states. The realizability of the bar states is enforced using the MCL technology for homogeneous hyperbolic systems. The forcing induced by stiff scattering is incorporated using Strang-type operator splitting. We use an explicit strong-stability-preserving Runge--Kutta method for the radiation transport subproblem and exact integration in the forcing steps, which guarantees the IDP property. The deposited dose is defined as the integral of a weighted zeroth moment over a bounded energy range. It is accumulated during the backward-in-energy evolution.
  Numerical experiments demonstrate that the proposed Strang-MCL method produces accurate and physically consistent dose distributions.
  
\end{abstract}

\begin{keyword}
  proton therapy, radiative transfer, realizable moment models,
  hyperbolic balance laws, finite elements,
  invariant domain preservation, flux limiting
\end{keyword}

\end{frontmatter}


\section{Introduction}

Proton therapy enables precise dose delivery to cancerous tissue while minimizing damage to surrounding healthy tissue. This remarkable capability arises from the so-called Bragg peak effect, a highly localized energy deposition occurring near the end of the proton range. Accurate prediction of the Bragg peak and the resulting dose distribution in patients is therefore essential for treatment planning. To date, Monte Carlo algorithms are widely regarded as the standard for clinical dose calculations~\cite{JANSON2024, Lin2021, Saini2018, verbeek2021}. Despite their high accuracy, the substantial computational effort required by these methods remains a major challenge for routine clinical use.

Deterministic models, such as the Boltzmann transport equation and its Fokker--Planck approximation for charged particles, offer a promising alternative for dose calculation in radiotherapy~\cite{bedford2019, gifford2006, Stammer2025, ulikema2012, vassiliev2010}. The physical processes governing dose distribution include energy loss due to ionization and the lateral spreading of proton beams caused by multiple Coulomb scattering. These effects are modeled by the stopping power~\cite{bortfeld1997, ulmer2007} and the scattering power~\cite{Gottschalk2009}, respectively. The proton fluence depends on space, energy, and direction of flight, resulting in a high-dimensional phase space.

Moment models, such as the $M_1$ model, provide a favorable compromise between accuracy and computational cost~\cite{berthon2011, duclous2010, pichard2016, frank2007}. The nonlinear $M_1$ system evolves the zeroth and first angular moments of the particle fluence. Since the number of unknowns exceeds the number of equations, the system must be closed by an additional relation that expresses the second angular moment in terms of the zeroth and first moments. Such closures can be constructed, e.g., using the maximum entropy principle~\cite{alldredge2012, brunner2000, brunner2001, coulombel2006, frank2007, frank2012, hauck2011, levermore1996, Levermore1984, minerbo1978, monreal2013, pichard2017}. Physical consistency requires that the reconstructed moments remain within a convex \emph{realizable set} of angular moments associated with nonnegative probability distributions~\cite{berthon2007, kershaw1976, olbrant2012}.  Entropy-based closures guarantee this property.

A discontinuous Galerkin (DG) method using slope limiters to ensure realizability for the time-dependent but energy-independent $M_1$ model can be found in \cite{chidyagwai2018}. The limiting strategies considered in \cite{chidyagwai2018} were found to introduce nonphysical disturbances. A realizable and nonoscillatory continuous finite element discretization of the same model was designed in our recent work~\cite{moujaes2026} using a combination of the monolithic convex limiting (MCL) framework~\cite{kuzmin2020} for hyperbolic systems of conservation laws and Patankar's method \cite{burchard2003} for positivity-preserving inclusion of source terms.

In the present work, we extend the steady-state MCL formulation to account for the energy dependence and the presence of the stopping power in the energy derivative. Treating the energy variable as a pseudo-time coordinate, we employ a backward marching method for energy discretization. Instead of treating the reactive forcing terms due to stiff scattering implicitly as in~\cite{moujaes2026}, we decouple them from the homogeneous $M_1$ system using
Strang-type operator splitting. The fractional-step Strang-MCL algorithm guarantees positivity preservation for the zeroth moment and enforces the \textit{realizable velocity} constraint for the first moment. These two moments can be interpreted as density and momentum, respectively. The application of MCL in the hyperbolic transport step and exact integration in the forcing step keep the nodal states in the realizable set of the $M_1$ model. From a formal mathematical perspective, our fully discrete scheme is provably \emph{invariant domain preserving} (IDP).

The deposited dose is computed by integrating the product of the zeroth moment and the stopping power over a finite energy interval. During the backward-in-energy evolution of the moments, the contribution of each energy step to the cumulative dose is calculated using the trapezoidal quadrature rule. The energy stepping is terminated at a small cutoff energy, at which point the remaining energy is assumed to be deposited locally. The results of our numerical studies for prototypical proton therapy scenarios are free of spurious oscillations and fully consistent with the underlying physics.

The remainder of this paper is organized as follows. Section~\ref{sec:modelling} introduces the high-dimensional proton transport model, the corresponding $M_1$ moment approximation, and an entropy-based closure. In Section~\ref{sec:loworder}, we present a low-order discretization that ensures the IDP property, i.e., realizability. The MCL algorithm that imposes global and local bounds on the quantities of interest in the flux-corrected high-order extension is described in Section~\ref{sec:MCL}. The procedure for dose calculation is outlined in Section~\ref{sec:dose}. Numerical results are reported in Section~\ref{sec:examples}, followed by conclusions in Section~\ref{sec:conclusion}.

\section{Proton transport modeling} 
\label{sec:modelling}

We begin with a review of two popular radiative transfer models for applications in proton therapy.

\subsection{Fokker-Planck equation}
\label{sec:physmodel}
The steady-state continuous slowing down Fokker--Planck approximation of the linear Boltzmann equation reads~\cite{bedford2019, frank2007, gifford2006}
\begin{equation}\label{eq:fp}
    \OO\cdot\nabla_\x \psi(\x, \OO, E) - \frac{\partial (S(\x, E)\psi(\x, \OO, E))}{\partial E} = \frac{T(\x, E)}{2}\Delta_\Omega \psi(\x, \OO, E), 
\end{equation}
where $\psi(\mathbf{x},\Omega,E)\geq 0$ denotes the proton fluence at position $\x\in\D\subset\R^d$, $d \in\{1,2,3\}$, moving in direction $\OO\in\Sp=\{\OO'\in\R^d: |\OO'|= 1\}$ with energy $E\in[0, E_{\max}]$. The physical meaning and modeling of the scaling functions $S(\x, E)$ and $T(\x, E)$ are explained below. For simplicity, we omit the explicit dependence on $\x$ in the models introduced below.

The stopping power $S(E):= -\frac{\diff E}{\diff x}$ represents the mean energy loss per unit path length.
There exist published data sets~\cite{berger1998} and parameterizations~\cite{newhauser2015} for $S(E)$. In this work, we adopt the Bragg--Kleeman rule~\cite{ashby2025, bortfeld1997, ulmer2007}, which approximates $S(E)$ by
\begin{equation}\label{eq:SBK}
    S(E) = \frac{1}{\beta p} E^{1-p}.
\end{equation}
This parametrization is derived from the continuous slowing-down range
\begin{equation}\label{eq:range}
    R(E) = \beta E^{p},
\end{equation}
which represents the range of a monoenergetic proton beam in a homogeneous medium. The parameters $\beta>0$ and $p\in[1,2]$ are typically fitted to experimental data and depend on the material~\cite{newhauser2015, ulmer2007}. For heterogeneous media, these parameters are commonly chosen to be piecewise constant in space, i.e., constant within each material slab~\cite{cox2024}. The Bragg--Kleeman model provides sufficient accuracy while remaining computationally simple compared to other parameterizations~\cite{ulmer2007}.

\begin{remark}
The numerical methods proposed below are compatible with any stopping power, provided it remains physically meaningful, i.e., strictly positive.
\end{remark}

The projected scattering power $T(E)>0$, also known as the angular diffusion coefficient, represents the rate at which a proton beam spreads laterally due to multiple small-angle Coulomb scatterings as it travels through a material.
We adopt the Rossi parameterization~\cite{Gottschalk2009}
\begin{equation}\label{eq:scatpow}
    T(E) = \left( \frac{E_s}{pv}\right)^2 \frac{1}{X_S},
\end{equation}
where $E_s= 15.0 \mev$, while $p$ and $v$ are the momentum and velocity of the proton, respectively.
The quantity $pv$, which depends on the proton energy $E$, is defined as
\begin{equation}\label{eq:pv}
    pv = \frac{\tau +2}{\tau +1}E,\qquad\tau =\frac{E}{mc^2},
\end{equation}
where $mc^2 = 938.272\,\mathrm{MeV}$ is the rest energy of protons.
The quantity $X_S$ denotes the scattering length~\cite{Gottschalk2009} and is given by
\begin{equation}\label{eq:scatlength}
    \frac{1}{\rho X_S} = \alpha N_A r_e^2 \frac{Z^2}{A}\left( 2\ln\left( 33219(AZ)^{-\frac13} \right)-1 \right),
\end{equation}
where $\alpha\approx \frac{1}{137}$ is the fine-structure constant, $N_A = 6.0221408\times 10^{23}$ is Avogadro's number, $r_e = 2.81796\times 10^{-13}\,\mathrm{cm}$ is the classical electron radius, and $\rho$, $A$, and $Z$ denote the mass density, atomic weight, and atomic number of the target material, respectively.
For compound or mixed materials, the scattering length can be computed using the Bragg rule~\cite{Gottschalk2009}
\begin{equation}\label{eq:braggrule}
    \frac{1}{\rho X_S} = \sum_m w_m \left( \frac{1}{\rho X_S} \right)_m,
\end{equation}
where $w_m$ is the fraction by weight of the $m$-th constituent, which we obtain from~\cite[Table 5.1]{IAEA.2024}.
Similarly to stopping power, the spatial dependence of the scattering power results from the material dependence of these quantities. Table~\ref{tab:BK_XS} summarizes the model parameters for selected materials.
\begin{table}[htbp]
\centering

\begin{tabular}{lcccc}
\hline
Material & $\beta$ & $p$ & $X_S$ & $\rho$ \\
\hline
Water   & 0.0022 & 1.77 & 46.88 & 1 \\
Muscle  & 0.0021 & 1.75 & 45.88 & 1.04 \\
Lung    & 0.0033 & 1.74 & 175.58 & 0.3\\
Bone    & 0.0011 & 1.77 & 17.93 &1.85\\
\hline
\end{tabular}
\caption{Bragg--Kleeman parameters~\cite[Table 1]{ashby2025}, scattering lengths obtained from~\eqref{eq:scatlength} and~\eqref{eq:braggrule}, and densities for selected materials.}
\label{tab:BK_XS}
\end{table}

The main quantity of interest for proton therapy is the dose~\cite{hensel2006,larsen1997}
\begin{equation}\label{eq:dose}
    D(\x) = \frac{1}{\rho(\x)}\int_0^{E_{\max}}\int_{\Sp} S(\x, E) \psi(\x, \OO, E)\,\diff\OO\,\diff E,
\end{equation}
which is absorbed by the medium.

\subsection{$M_1$ moment model}
\label{sec:momentmodel}

Numerical solution of the Fokker–Planck equation~\eqref{eq:fp} is computationally expensive due to the high dimensionality of the domain $\D\times \Sp \times [0,E_{\max}]$. The $M_N$ model approximates \eqref{eq:fp} by a nonlinear system of equations for the first $N+1$ angular moments of $\psi$ defined by
    \begin{equation}\label{eq:psin}
    \psi^{(n)} = \int_{\Sp}\underbrace{\OO \otimes \cdots \otimes \OO}_{n\text{ times}}\psi(\OO)\,\diff\OO,\qquad n = 0,\ldots,N.
    \end{equation}
The corresponding balance laws are derived by taking angular moments of
\eqref{eq:fp}.
A suitable closure is required to approximate
$\psi^{(N+1)}$ in terms of $\psi^{(0)},\ldots,\psi^{(N)}$.
Hereafter, boldface notation is used for angular moments of degree $n>0$,
which represent vector or tensor fields.

Focusing on the case $N=1$, we consider
the $M_1$ system~\cite{duclous2010}
\begin{align}
        \nabla \cdot \po - \frac{\partial\left(S\pz\right)}{\partial E} &= 0,\label{eq:M1_0moment}\\
        \nabla \cdot \pt - \frac{\partial \left(S\po\right)}{\partial E} &= -T\po,\label{eq:M1_1moment}
\end{align}
which consists of coupled equations for the zeroth and first moments
\begin{align*}
    \pz = \pz(\x,E) &= \int_{\Sp} \psi(\x, \OO, E)\,\diff\OO\,\in\R,\\
    \po = \po(\x,E) &= \int_{\Sp} \OO \psi(\x, \OO, E)\,\diff\OO\,\in\R^d. 
\end{align*}
The vanishing right-hand side of~\eqref{eq:M1_0moment} follows from the self-adjointness of the Laplace--Beltrami operator $\Delta_{\OO}$ on the unit sphere. Furthermore, each component of the direction vector $\OO=(O_1,\ldots,O_d)$ can be expressed as a linear combination of spherical harmonics of degree $\ell= 1$. These harmonics are eigenfunctions of $\Delta_{\OO}$ with eigenvalue $\lambda = -\ell(\ell +1) = -2$. The reactive source term of equation~\eqref{eq:M1_1moment} arises from the application of $\Delta_{\OO}$ to the components of $\OO$.

The system of angular moment equations \eqref{eq:M1_0moment} and \eqref{eq:M1_1moment}
can be cast in the form 
\begin{equation}\label{eq:M1}
    \dv \f(u) - \frac{\partial\left(Su\right)}{\partial E} = -\sigma u,
\end{equation}
where the vector of unknowns $u$ and the matrix $\f(u)$ of corresponding fluxes are given by
\begin{equation*}
  u=\left(\begin{array}{cc}
         \pz\\
         \po
  \end{array}\right)\in\R^{d+1},
  \qquad
    \f(u) = \left(\begin{array}{cc}
         \po\\
         \pt
    \end{array}\right)\in\R^{d\times(d+1)},
\end{equation*}
and $\sigma = \mathrm{diag}(0, T, \ldots, T) \in\R^{(d+1)\times(d+1)}$.

Note that the $n$-th moment is transported via the $(n+1)$-st moment. 
Thus, we require a closure relation for $\pt = \pt(\pz,\po)$ to approximate the  second moment 
\begin{equation*}
    \pt = \pt(\x,E) = \int_{\Sp} \OO\otimes\OO  \psi(\x, \OO, E)\,\diff\OO\,\in\R^{d\times d}.
\end{equation*}

To ensure the physical validity of the $M_1$ model, the closure must guarantee that if $\pz$ and $\po$ are moments of a nonnegative angular distribution $\psi\geq 0$, the reconstruction $\pt$ corresponds to the second moment of~$\psi$. Such a closure is called realizable in the radiative transfer literature. A common approach for recovering $\pt$ is via (approximate) maximum entropy reconstruction~\cite{alldredge2012, chidyagwai2018, coulombel2006, Levermore1984, pichard2017}. To avoid solving a potentially ill-conditioned optimization problem at each point $(\x,E)\in\D\times[0,E_{\max}]$, we employ the widely used realizable approximation~\cite{Levermore1984}
\begin{equation}\label{eq:M1closure}
    \pt = \mathbf{D}\left(\vel\right)\pz, \quad \vel = \frac{\po}{\pz},
\end{equation}
where
\begin{equation}\label{eq:D}
    \mathbf{D}(\vel) = \frac{1-\chi(|\vel|)}{2} I_d + \frac{3\chi(|\vel|)-1}{2} \frac{\vel \otimes \vel}{|\vel|^2}
\end{equation}
is the Eddington tensor and
\begin{equation}\label{eq:chi}
    \chi(f) = \frac{3+4f^2}{5+2\sqrt{4-3f^2}}
\end{equation}
is the Eddington factor.

It is easy to verify that the zeroth and first moments of a nonnegative angular function $\psi$ satisfy
\begin{equation}\label{eq:m01realizable}
    \pz \geq 0\quad\text{and}\quad|\po|\leq \pz,
\end{equation}
respectively, while $\pz = 0$ if and only if $\psi(\OO) \equiv 0$. 
In this trivial case, $|\po|=\pz$. Equality $|\po| = \pz$ occurs only for a perfectly collimated beam, corresponding to an angular delta distribution $\psi(\OO) = \delta(\OO-\OO')$ for some $\OO' \in \Sp$~\cite{kershaw1976}.
For $|\po| = \pz$, the directional Jacobian of the flux function, $\f'_{\n}(u) = \frac{\partial}{\partial u} (\f(u)\cdot \n)$, becomes non-diagonalizable, and the $M_1$ system~\eqref{eq:M1} is no longer hyperbolic~\cite{chidyagwai2018}.
Accordingly, we define the \textit{realizable set} of physically admissible states as
\begin{equation}\label{eq:R1}
    \mathcal{R}_1 = \left\{ (\pz, \po)^{\top} \in \R^{d+1} : \pz > 0, \ |\po| < \pz \right\},
\end{equation}
which corresponds to the set of moments of nonnegative, nontrivial $L^2(\Sp)$ functions.
Note that the realizable set is a convex cone. 
Furthermore, we refer to moments belonging to $\mathcal{R}_1$ as realizable.

Let $u= (\pz,\po)^\top\in\mathcal{R}_1$ and reconstruct $\pt(u)$ using \eqref{eq:M1closure}--\eqref{eq:chi} with
\begin{equation*}
    f^2\leq\chi(f)\leq 1\quad \text{for }f\in[0,1).
\end{equation*}
Under these conditions, $\pz,\po,$ and $\pt$ represent the zeroth, first, and second moments of a nonnegative angular distribution, respectively~\cite{Levermore1984}.

In view of~\eqref{eq:dose}, the absorbed dose can be written in terms of the zeroth moment as follows \cite{duclous2010}: 
\begin{equation}\label{eq:M1dose}
    D(\x) = \frac{1}{\rho(\x)}\int_0^{E_{\max}} S(\x, E) \pz(\x, E)\,\diff E.
\end{equation}

\section{Discretization and methodology}
To construct a realizability-preserving continuous finite element discretization
of \eqref{eq:M1}, we first derive a low-order scheme that is invariant-domain preserving (IDP), in the sense that the approximate nodal states remain in the realizable set $\mathcal{R}_1$. We then incorporate high-order correction terms and use the monolithic convex limiting (MCL) methodology~\cite{kuzmin2020} to enforce the IDP property. For the time-dependent and energy-independent $M_1$ model, such algorithms were designed in our previous work~\cite{moujaes2026}. In this section, we adapt them to the structure of system~\eqref{eq:M1}. The product $Su$ is evolved backward in energy using Strang splitting to decouple the forcing terms arising from scattering. Exact energy integration for the source-term subproblems provides a more accurate IDP treatment than the algorithm employed in~\cite{moujaes2026}. An explicit strongly-stability-preserving Runge--Kutta (SSP-RK) method is used for the homogeneous subsystem of~\eqref{eq:M1}. The dose~\eqref{eq:M1dose} is decomposed into integrals over individual energy evolution steps. These integrals are approximated using the trapezoidal rule.

\subsection{High-order method}\label{sec:highorder}

Let a boundary condition of the form $\f(u)\cdot \n=\mathcal{F}(u, \hat u; \n)$ be imposed weakly on the boundary $\Gamma = \partial \D$ of the spatial domain $\D\subset\R^3$. Choosing a test function $w$, we construct the weak form
\begin{equation}\label{eq:weakform}
    \int_\D w\left(  \dv \f(u) - \frac{\partial(Su)}{\partial E} + \sigma u \right)\dx = \int_\Gamma w\left(\f(u)\cdot \n - \mathcal{F}(u, \hat u; \n) \right)\,\diff\s
\end{equation}
 of the $M_1$ system~\eqref{eq:M1}. The boundary term is defined using the
 global Lax--Friedrichs (GLF) flux 
\begin{equation*}
    \mathcal{F}(u, \hat u; \n) = \frac{\f(u) + \f(\hat u)}{2} \cdot \n -\frac{\lambda_{\max}}{2} (\hat u - u),
\end{equation*} 
where $\lambda_{\max}=1$ is a global upper bound for the maximum wave speed~\cite{berthon2007, chidyagwai2018, olbrant2012}. The external state $\hat u$ of the boundary condition
corresponds to a proton beam. 

We discretize the weak form~\eqref{eq:weakform} in space using a conforming mesh $\mathcal{T}_h =\{K_e\}_{e=1}^{E_h}$ with cells $K_1,\ldots,K_{E_h}$ and vertices $\x_1,\ldots,\x_{N_h}$. The Lagrange basis functions of a globally continuous, piecewise linear ($\mathbb{P}_1$) or multilinear ($\mathbb{Q}_1$) finite element approximation are denoted by $\varphi_1,\ldots,\varphi_{N_h}$. They possess the interpolatory property $\varphi_i(\x_j) =\delta_{ij}$ and span the space $V_h\subset H^1(\mathcal D)\cap C(\bar{\mathcal D})$.

A \emph{group finite element} approximation (cf. \cite{barrenechea2017b,fletcher1983})
to the conserved product $Su$ is defined by
\begin{equation}\label{eq:Suh}
   (Su)_h(\x,E) = \sum_{j = 1}^{N_h} (Su)_j(E)\varphi_j(\x).
\end{equation}
In a similar vein, the nonlinear flux function $\f(u(\x,E))$ of the $M_1$ model is approximated by
\begin{equation}\label{eq:GFE}
    \f_h(\x,E)= \sum_{j = 1}^{N_h} \f_j(E)\varphi_j(\x),\quad \f_j(E) = \f(u_j(E)).
\end{equation}
The coefficients $(Su)_j(E) = (Su)_h(\x_j, E)$ of $(Su)_h$ are evolved directly, whereas
the nodal states $u_j(E) = u_h(\x_j, E)$ of the corresponding moment approximation
\begin{equation}\label{eq:uh}
    u_h(\x,E) = \sum_{j=1}^{N_h} u_j(E)\varphi_j(\x),\qquad
 u_j(E) = \frac{(Su)_j(E)}{S_j(E)},\quad S_j(E) = S(\x_j, E) 
\end{equation}
are recovered from the main discrete unknowns $(Su)_j(E)$  via division by the nodal stopping power $S_j(E)$. Since the realizable set $\mathcal{R}_1$ is a convex cone, $(Su)_j\in\mathcal{R}_1$ if and only if $u_j\in\mathcal{R}_1$.

We introduce the index sets $\N_i = \{j\in \{1,..., N_h\}: \mathrm{supp}(\varphi_j)\cap \mathrm{supp}(\varphi_i)\not =\emptyset\}$ and $\N_i^* = \N_i\setminus\{i\}$ to define the computational stencils of node $i$.
Using a basis function $w = \varphi_i$ as test function and substituting the finite element approximations \eqref{eq:Suh}--\eqref{eq:GFE} into~\eqref{eq:weakform}, we obtain 
\begin{equation}\label{eq:CGsemidiscrete}
    -\sum_{j\in \N_i} m_{ij} \frac{\diff (Su)_j}{\diff E} = b_i(u_h, \hat u) - \sum_{j\in \N_i} [\f_j \cdot \con_{ij} + m_{ij}^{\sigma}u_j].
\end{equation}
The coefficients of this semi-discrete backward-in-energy evolution equation are given by
\begin{align*}
    m_{ij} &=\int_\D \varphi_i \varphi_j\,\diff\x,\quad \con_{ij} = \int_\D \varphi_i \nabla\varphi_j\,\diff\x,\\
    m_{ij}^{\sigma} &=\mathrm{diag}(0, m_{ij}^T,\ldots,m_{ij}^T),\quad  m_{ij}^T= \int_\D T \varphi_i \varphi_j\,\diff\x. 
\end{align*}
The weakly imposed boundary condition is taken into account via
\begin{equation*}
     b_i(u_h, \hat u) = \int_\Gamma \varphi_i \left(\f(u_h) \cdot \n - \mathcal{F}(u_h, \hat u;\n)\right)\,\diff\s.
\end{equation*}
Since the standard Galerkin discretization \eqref{eq:CGsemidiscrete} is generally not
IDP, we will modify it using a customized version of
the MCL procedures developed in \cite{kuzmin2020,moujaes2026}.

\subsection{Low-order method}\label{sec:loworder}

A fundamental building block of the MCL algorithm to be designed is a low-order scheme that
provides provable IDP properties. We derive it from \eqref{eq:CGsemidiscrete}
using \emph{mass lumping} (inexact nodal quadrature) and artificial \emph{graph viscosity}
of GLF type. The lumped-mass approximations
\begin{equation*}
    \sum_{j\in \N_i} m_{ij} \frac{\diff (Su)_j}{\diff E} \approx m_i \frac{\diff (Su)_i}{\diff E}, \qquad \sum_{j\in \N_i} m_{ij}^{\sigma} u_j \approx m_i^{\sigma}u_i
\end{equation*}
are defined using the row sums of the corresponding consistent mass matrices, i.e.,
\begin{align*}
    m_i &= \sum_{j\in\N_i} m_{ij} =  \int_{\D}\varphi_i\,\diff\x,\qquad
    m_i^{\sigma} = \mathrm{diag}(0, m_{i}^T,\ldots,m_{i}^T), \qquad m_i^T= \sum_{j\in\N_i} m_{ij}^{T} =  \int_{\D}T\varphi_i\,\diff\x.
\end{align*}
The boundary term $b_i(u_h, \hat u)$ of problem \eqref{eq:CGsemidiscrete}
is approximated by its lumped counterpart
\begin{equation}\label{eq:lumpedbdr}
   \tilde{b}_i(u_i, \hat u_i) = \int_\Gamma \varphi_i\left(\f_i \cdot \n - \mathcal{F}(u_i, \hat u_i;\n)\right)\,\diff\s. 
\end{equation}
Finally, low-order stabilization via dissipative numerical fluxes of the form $d_{ij}(u_j -u_i)$
is incorporated into the semi-discrete finite element GLF scheme
\begin{equation}\label{eq:LO}
    -m_i \frac{\diff (Su)_i}{\diff E} = \tilde{b}_i(u_i, \hat u_i)+ \sum_{j\in \N_i^*}\left[d_{ij} (u_j -u_i) - (\f_j-\f_i)\cdot \con_{ij}\right] -  m_i^{\sigma}u_i.
\end{equation}
The GLF graph viscosity coefficients~\cite{kuzmin2010a,moujaes2026}
\begin{equation*}
    d_{ij} = \begin{cases}
        \lambda_{\max}\max\{|\con_{ij}|, |\con_{ji}|\} &\text{if }j\in\N_i^*,\\
        -\sum_{k\in\N_i^*}d_{ik} &\text{if } j=i,\\
        0&\text{otherwise}
    \end{cases}
\end{equation*}
are defined using the global bound $\lambda_{\max} = 1$ for the realizable
maximum speed of the $M_1$ model.

To show the IDP property of a fully discrete version of~\eqref{eq:LO}, we introduce the \textit{bar states}~\cite{guermond2016, kuzmin2020, kuzmin2023}
\begin{equation}\label{eq:LObarstates}
    \bar u_{ij} = \frac{u_i + u_j}{2} - \frac{(\f_j-\f_i)\cdot \con_{ij}}{2d_{ij}}.
\end{equation}
As shown in~\cite{moujaes2026}, these intermediate states
 belong to $\mathcal{R}_1$ if $u_i, u_j\in\mathcal{R}_1$ and $d_{ij}\geq |\con_{ij}|$.
 For simplicity, we assume $i$ to be an interior node
 and use~\eqref{eq:LObarstates} to write~\eqref{eq:LO} with $\tilde b_i(u_i,\hat u_i) = 0$
 in the form 
\begin{equation}\label{eq:LObsform}
     -m_i \frac{\diff (Su)_i}{\diff E} = \sum_{j\in\N_i^*}[ 2d_{ij} (\bar{u}_{ij} -u_i)] - m_i^{\sigma}u_i.
\end{equation}
It is convenient to 
decompose the bar state form~\eqref{eq:LObsform} into the following two subproblems:
\begin{align}
    -m_i \frac{\diff (Su)_i}{\diff E} &=- m_i^{\sigma}u_i,\label{eq:scatteringprob}\\
    -m_i \frac{\diff (Su)_i}{\diff E} &= \sum_{j\in\N_i^*}[ 2d_{ij} (\bar{u}_{ij}-u_i)] \label{eq:transportprob}.
\end{align}
We use this decomposition and  the
symmetric Strang splitting procedure to advance the numerical solution from the energy level $E^{n+1}$ to $E^{n} = E^{n+1} - \Delta E$ in three steps~\cite[Sec. 6.2.7]{kuzmin2014a}:

\begin{enumerate}
\item \label{it:scat}
Solve the scattering subproblem~\eqref{eq:scatteringprob} over a half-step in energy from $E^{n+1}$ to
$E^{n+\frac12}=E^{n+1}-\frac{\Delta E}{2}$.

\item \label{it:transport}
Starting from the intermediate solution obtained in Step~\ref{it:scat}, solve the homogeneous radiation transport subproblem~\eqref{eq:transportprob} over a full step in energy from $E^{n+1}$ to $E^{n}$.

\item
Starting from the intermediate solution obtained in Step~\ref{it:transport},
solve the scattering subproblem~\eqref{eq:scatteringprob} over the second half-step in energy from $E^{n+\frac12}$ to the final energy $E^{n}$.
\end{enumerate}

We first consider the scattering subproblem~\eqref{eq:scatteringprob}, which is linear and decouples into independent problems at each node.
Let the nodal state
\begin{equation*}
    (Su)_i^{n+1} = ((S\pz)_i^{n+1},(S\po)_i^{n+1})^\top
\end{equation*} 
at energy $E^{n+1}$ be realizable.
Since $m_i^\sigma = \mathrm{diag}(0,m_i^T, \ldots ,m_i^T)$, the equation for the zeroth moment has a vanishing right-hand side, and thus $(S\pz)_i^{n+1/2} = (S\pz)_i^{n+1}$. The evolution
equation
\begin{equation}\label{eq:scatequation}
    -m_i \frac{\diff (S\po)_i}{\diff E} =-\frac{m_i^{T}}{S_i} (S\po)_i
\end{equation}
for the first moment can be integrated exactly over $[E^{n+1/2}, E^{n+1}]$ and rearranged to obtain
\begin{equation}\label{eq:scatteringstep}
    (S\po)_i^{n+1/2} = \exp\left( -\frac{1}{m_i}\int_{E^{n+1/2}}^{E^{n+1}} \frac{m_i^T}{S_i}\,\diff E\right) (S\po)_i^{n+1}.
\end{equation}
In practice, the energy integral is approximated using the midpoint rule.
This scattering update preserves realizability because  $(S\pz)_i^{n+1/2} = (S\pz)_i^{n+1} >0 $ and 
\begin{equation*}
    \left|(S\po)_i^{n+1/2}\right|\leq \left|(S\po)_i^{n+1}\right| < (S\pz)_i^{n+1} = (S\pz)_i^{n+1/2} .
\end{equation*}

Next, we consider the transport subproblem~\eqref{eq:transportprob} on the energy interval from $E^{n+1}$ to $E^n$ starting from the intermediate solution $(Su)_i^{n+1/2}\in\mathcal{R}_1$ obtained with~\eqref{eq:scatteringstep}.
We discretize~\eqref{eq:transportprob} using an explicit SSP-RK method backwards.
Each backward Euler stage can be written as
\begin{equation}\label{feupdate}
\begin{split}
    (Su)_i^{\mathrm{SSP}} =& (Su)_i + \frac{\Delta E}{m_i}\sum_{j\in\N_i^*}[ 2d_{ij} (\bar{u}_{ij} -u_i)]\\
    =&  S_i\left[ \left(1-\frac{\Delta E}{m_i  S_i} \sum_{j\in\N_i^*} 2d_{ij}\right) u_i +\frac{\Delta E}{m_i  S_i} \sum_{j\in\N_i^*} 2d_{ij} \bar{u}_{ij}\right],
\end{split}
\end{equation}
where $(Su)_i$ is defined at $E^{n+1}$ and $(Su)_i^{\mathrm{SSP}}$ at $E^n$. 
Under the CFL-like condition 
\begin{equation}\label{eq:CFL}
    \frac{2\Delta E}{m_i S_i}\sum_{j\in\N_i^*} d_{ij}\leq 1,
\end{equation}
the result $(Su)_i^{\mathrm{SSP}}$ of the explicit update \eqref{feupdate} is realizable, because it represents  a scaled convex combination of the realizable states $u_i$ and $\bar{u}_{ij}$, $j\in\N_i^*$. 

\begin{remark}
 Similarly to the right-hand side of~\eqref{eq:transportprob}, 
 a nonvanishing boundary term~\eqref{eq:lumpedbdr} 
 can be written in an IDP bar state form~\cite{hajduk2021, kuzmin2023, moujaes2025}.
 The realizable external states $u_j=\hat u_i$ with indices $j>N_h$ represent the Riemann data of the
 weakly imposed boundary condition.
 The corresponding generalization of~\eqref{eq:transportprob} preserves the realizability of nodal moment states under a suitable restriction on $\Delta E$.
\end{remark}


Since the solvers for individual subproblems of the low-order method using Strang splitting
are IDP, realizability is also guaranteed for the final solution at the energy level $E^n$. 

\subsection{Monolithic convex limiting}\label{sec:MCL}
The high-order spatial semi-discretization~\eqref{eq:CGsemidiscrete}
can be recovered from the low-order method~\eqref{eq:LO} by adding  antidiffusive fluxes
$f_{ij}$ that correct the mass lumping error and offset the diffusive fluxes
$d_{ij} (u_j - u_i)$. Let $\mathrm{d}_E(Su)_i = \frac{\diff (Su)_i}{\diff E}$ denote
the nodal energy derivatives corresponding to~\eqref{eq:CGsemidiscrete}. Then
\begin{equation}\label{eq:rawADF}
  f_{ij} = - m_{ij}\left( \mathrm{d}_E(Su)_i -
  \mathrm{d}_E{(Su)}_j\right) + (d_{ij} + m_{ij}^{\sigma})( u_i - u_j). 
\end{equation}
To avoid solving a linear system with the consistent mass matrix $(m_{ij})_{i,j = 1}^{N_h}$ and to incorporate high-order stabilization into $f_{ij}$, we approximate  $\mathrm{d}_E(Su)_i$ by
(cf. \cite{kuzmin2020,kuzmin2023,lohmann2019})
\begin{equation*}
  \mathrm{d}_E{(Su)}_i^L = \frac{1}{m_i}\left(m_i^{\sigma} u_i-
  \sum_{j\in \N_i^*}\left[d_{ij} (u_j -u_i) - (\f_j-\f_i)\cdot \con_{ij}\right]\right).
\end{equation*}
Substituting the antidiffusive fluxes~\eqref{eq:rawADF} into~\eqref{eq:LO} recovers the stabilized high-order target scheme
\begin{equation}\label{eq:HOtarget}
\begin{split}
     -m_i \frac{\diff (Su)_i}{\diff E} &= \sum_{j\in \N_i^*}\left[d_{ij} (u_j -u_i) - (\f_j-\f_i)\cdot \con_{ij} + f_{ij}\right] - m_i^{\sigma} u_i\\
      &= \sum_{j\in\N_i^*}[ 2d_{ij} (\bar{u}^H_{ij} -u_i)] - m_i^{\sigma} u_i.
\end{split}
\end{equation}
The high-order bar states
\begin{equation}\label{eq:HObs}
    \bar{u}^H_{ij} = \bar{u}_{ij} + \frac{f_{ij}}{2d_{ij}}
\end{equation}
generally do not belong to $\mathcal{R}_1$.
The monolithic convex limiting (MCL) strategy proposed in~\cite{kuzmin2020} replaces the raw antidiffusive fluxes $f_{ij} = - f_{ji}$ by their limited counterparts $f_{ij}^* = -f_{ji}^*$ such that the realizability of $\bar{u}_{ij},\bar{u}_{ji}\in\mathcal{R}_1$ is preserved by the
flux-corrected bar states
\begin{equation}\label{eq:limitedBS}
    \bar{u}^*_{ij} = \bar{u}_{ij} + \frac{f_{ij}^*}{2d_{ij}},\qquad \bar{u}^*_{ji} = \bar{u}_{ji} - \frac{f_{ij}^*}{2d_{ij}}.
\end{equation}

In the MCL version, the radiation transport subproblem of the Strang splitting method is discretized
using~\eqref{eq:transportprob} with $\bar{u}_{ij}$ replaced by $\bar{u}^*_{ij}$. Since the structure of the low-order scheme is preserved, the backward Euler stages of the explicit SSP-RK method are IDP under the CFL-like condition~\eqref{eq:CFL}.

In addition to enforcing the physical admissibility conditions $\bar{u}^*_{ij},\bar{u}^*_{ji}
\in\mathcal{R}_1$, a well designed flux limiter should effectively suppress spurious
oscillations in the neighborhood of shocks and steep gradients. The MCL procedure
employed in \cite{moujaes2025} achieves numerical admissibility by imposing local
discrete maximum principles on the scalar-valued components of the bar states~\eqref{eq:limitedBS}.

Let $\phi_i$ be a component of $u_i = (\pz_i, \psi_{i,1}^{(1)},\,\ldots,\,\psi_{i,d}^{(1)})^\top\in\mathcal{R}_1$. The corresponding low-order bar states and raw antidiffusive fluxes
are denoted by $\bar\phi_{ij}$ 
and $f_{ij}^\phi$, $j\in\N_i^*$, respectively. We formulate the numerical admissibility conditions
\begin{equation}\label{eq:localbp}
  \begin{split}
    \phi_i^{\min}\leq\bar\phi_{ij}^* &= \bar\phi_{ij} +\frac{f_{ij}^{\phi,*}}{2d_{ij}} \leq \phi_i^{\max},\\
      \phi_j^{\min}\leq\bar\phi_{ji}^* &= \bar\phi_{ji} -\frac{f_{ij}^{\phi,*}}{2d_{ij}} \leq \phi_j^{\max}
    \end{split}
\end{equation}
for the limited antidiffusive fluxes $f_{ij}^{\phi,*}$ using the local bounds
\begin{equation}\label{eq:localbounds}
    \phi_i^{\max} = \max\left\{ \max_{j\in\N_i} \phi_j, \max_{j\in\N_i^*} \bar\phi_{ij} \right\},\quad
    \phi_i^{\min} = \min\left\{ \min_{j\in\N_i} \phi_j, \min_{j\in\N_i^*} \bar\phi_{ij} \right\}.
\end{equation}
The inequality constraints \eqref{eq:localbp} can be rearranged to 
\begin{equation}\label{eq:fluxbounds}
  \begin{split}
    2d_{ij}\left(\phi_i^{\min} -\bar\phi_{ij}\right)\leq f_{ij}^{\phi,*}&\leq 2d_{ij}\left( \phi_i^{\max} -\bar\phi_{ij}\right),\\
    2d_{ij}\left(\phi_j^{\min}-\bar\phi_{ji}
    \right)\leq -f_{ij}^{\phi,*}&\leq 2d_{ij}
    \left(\phi_j^{\max}-\bar\phi_{ji} \right).
    \end{split}
\end{equation}
Introducing the \emph{bounding fluxes}
\begin{align}
    f_{ij}^{\phi,\min} &= -f_{ji}^{\phi,\max} =  2d_{ij}\max \left\{\phi_i^{\min} -\bar\phi_{ij}, \bar\phi_{ji} -\phi_j^{\max}\right\}\leq 0,\\
    f_{ij}^{\phi,\max} &= -f_{ji}^{\phi,\min} = 2d_{ij}\min \left\{\phi_i^{\max} -\bar\phi_{ij}, \bar\phi_{ji} -\phi_j^{\min}\right\}\geq 0,
\end{align}
we set
\begin{equation}\label{eq:limitedflux}
    f_{ij}^{\phi,*} = \max\{f_{ij}^{\phi,\min}, \min\{f_{ij}^{\phi,\max}, f_{ij}^{\phi}\}\}.
\end{equation}
This adjustment of $f_{ij}^{\phi}$
ensures the discrete conservation property $f_{ij}^{\phi,*} = -f_{ji}^{\phi,*}$ and
the validity of \eqref{eq:localbp} for the bar states. The positivity of the zeroth moment is guaranteed too,
since $\bar\psi^{(0),*}_{ij}\geq\psi_i^{(0),\min}>0$ by construction of the local bounds
\eqref{eq:localbounds}.
However, enforcing the realizable velocity constraints
\begin{equation*}
    \left| \bar{\boldsymbol{\psi}}^{(1),*}_{ij} \right| \leq \bar\psi^{(0),*}_{ij},\qquad 
    \left| \bar{\boldsymbol{\psi}}^{(1),*}_{ji} \right| \leq \bar\psi^{(0),*}_{ji}
\end{equation*}
requires an additional flux limiting step, which was originally proposed in~\cite{moujaes2026} and inspired by a positivity fix for the specific internal energy of the compressible Euler equations~\cite{kuzmin2020,kuzmin2023}.

We denote by $f_{ij}^* = (f_{ij}^{*(0)}, \f_{ij}^{*(1)})^\top$ the limited antidiffusive fluxes whose individual components are given by~\eqref{eq:limitedflux}. To ensure that the maximum speed corresponding to the final bar state
\begin{equation}\label{eq:HOIDPbarstats}
    \bar{u}_{ij}^{\mathrm{IDP}} = \bar{u}_{ij} +\frac{\alpha_{ij}^{\mathrm{IDP}} f_{ij}^*}{2d_{ij}}
\end{equation}
remains bounded by $\lambda_{\max}=1$, we apply a correction factor
$\alpha_{ij}^{\mathrm{IDP}}\in[0,1]$ such that $\alpha_{ij}=\alpha_{ji}$ and
$\bar{u}_{ij}^{\mathrm{IDP}},\bar{u}_{ji}^{\mathrm{IDP}}\in\mathcal{R}_1$. Written in terms of moments, the IDO constraint for $\bar{u}_{ij}^{\mathrm{IDP}}$ becomes
\begin{equation*}
    \left|\bar{\bm\psi}^{(1)}_{ij} + \frac{\alpha_{ij}^{\mathrm{IDP}} \f_{ij}^{*(1)}}{2d_{ij}}\right|^2 
    < \left(\bar\psi^{(0)}_{ij} + \frac{\alpha_{ij}^{\mathrm{IDP}} f_{ij}^{*(0)}}{2d_{ij}}\right)^2.
\end{equation*}
This is a quadratic inequality constraint of the form
\begin{equation}\label{eq:pij<qij}
    P_{ij}(\alpha_{ij}^{\mathrm{IDP}}) < Q_{ij},
\end{equation}
where 
\begin{equation*}
    P_{ij}(\alpha) = \left(\left|\f_{ij}^{*(1)}\right|^2 - \left(f_{ij}^{*(0)}\right)^2 \right)\alpha^2 
    + 4 d_{ij} \left(\bar{\bm\psi}^{(1)}_{ij}\cdot \f_{ij}^{*(1)} -  \bar\psi^{(0)}_{ij}f_{ij}^{*(0)}\right)\alpha,
\end{equation*}
\begin{equation*}
    Q_{ij} = (2d_{ij})^2\left(
    \left(\bar\psi^{(0)}_{ij}\right)^2 - \left|\bar{\bm\psi}^{(1)}_{ij}\right|^2\right) > 0.
\end{equation*}
For any $\alpha\in[0,1]$, the estimate $\alpha^2\leq \alpha$  implies
 $P_{ij}\leq \alpha R_{ij}$, where
\begin{equation*}
    R_{ij} = \max\left\{  0, \left|\f_{ij}^{*(1)}\right|^2 - \left(f_{ij}^{*(0)}\right)^2  \right\} + 4 d_{ij} \left(\bar{\bm\psi}^{(1)}_{ij}\cdot \f_{ij}^{*(1)} -  \bar\psi^{(0)}_{ij}f_{ij}^{*(0)}\right).
\end{equation*}
To enforce the strict inequality~\eqref{eq:pij<qij}, we replace $\tilde Q_{ij}$
by $\tilde Q_{ij} = (1-\varepsilon)Q_{ij}>0$ with
$\varepsilon = 10^{-15}$. The IDP correction factor
\begin{equation*}
    \alpha_{ij}^{\mathrm{IDP}} = \begin{cases}
        \min\left\{ \frac{\tilde Q_{ij}}{R_{ij}} , \frac{\tilde Q_{ji}}{R_{ji}} \right\} &\text{if } R_{ij} > \tilde Q_{ij}, R_{ji} >\tilde Q_{ji},\\
        \frac{\tilde Q_{ij}}{R_{ij}} &\text{if }R_{ij} > \tilde Q_{ij}, R_{ji} \leq  \tilde Q_{ji},\\
        \frac{\tilde Q_{ji}}{R_{ji}} &\text{if }R_{ij} \leq  \tilde Q_{ij}, R_{ji} >  \tilde Q_{ji}.\\
        1 &\text{otherwise}
    \end{cases}
\end{equation*}
satisfies $\alpha_{ij}^{\mathrm{IDP}}=\alpha_{ji}^{\mathrm{IDP}}$. Furthermore,
 the physical admissibility conditions
\begin{equation*}
    P_{ij}(\alpha_{ij}^{\mathrm{IDP}}) \leq \alpha_{ij}^{\mathrm{IDP}} R_{ij} \leq \tilde Q_{ij} < Q_{ij}, \qquad P_{ji}(\alpha_{ij}^{\mathrm{IDP}}) \leq \alpha_{ij}^{\mathrm{IDP}} R_{ji} \leq \tilde Q_{ji}< Q_{ji}
\end{equation*}
hold for the final bar states~\eqref{eq:HOIDPbarstats}.
Hence, substituting $\bar{u}_{ij}^{\mathrm{IDP}}$ for $\bar u_{ij}$ in~\eqref{eq:transportprob} yields a numerically admissible and realizability preserving high-order extension
\begin{equation}\label{eq:sdMCL}
\begin{split}
     -m_i \frac{\diff (Su)_i}{\diff E} &= \sum_{j\in \N_i^*}\left[d_{ij} (u_j -u_i) - (\f_j-\f_i)\cdot \con_{ij} + \alpha_{ij}^{\mathrm{IDP}}f_{ij}^*\right]\\
      &= \sum_{j\in\N_i^*}[ 2d_{ij} (\bar{u}^{\mathrm{IDP}}_{ij} -u_i)]
\end{split}
\end{equation}
of the low-order semi-discrete scheme~\eqref{eq:transportprob}
for the transport subproblem of the Strang splitting algorithm.
The energy stepping methods remain unchanged for all subproblems.

\subsection{Dose calculation}\label{sec:dose}
The proposed method treats energy as a pseudo-time variable and advances the numerical solution backward in energy until the threshold $E_{\min} = 10^{-5}\mev$ is reached.
The numerical solution and the contributions to the dose are computed at the discrete energy levels
\begin{equation*}
    E_{\max} = E^N > E^{N-1}>\cdots>E^1 = E_{\min}.
\end{equation*}
Stopping at a strictly positive minimum energy $E_{\min}$ avoids the breakdown of physical models, e.g., for the stopping power~\eqref{eq:SBK} and the scattering power~\eqref{eq:scatpow}.
This ensures that no division by $\infty$ occurs when computing the integrals in the scattering step~\eqref{eq:scatteringstep} or the nodal states $(Su)_i/S_i$ in the transport step.
The choice of $E_{\max}$ depends on the boundary conditions.

In a typical numerical experiment for a box domain $\mathcal D\subset\R^d$ 
with inflow boundary $$\Gamma_{\mathrm{in}} = \{\x =(x_1,\ldots,x_d)\in\Gamma: x_1=0\},$$
we prescribe the zeroth moment 
\begin{equation}\label{eq:bdrcon0}
\hat \psi^{(0)}(\x, E) = \psi_0 \;
\frac{1}{\sqrt{2\pi}\, \sigma_E} \exp\Big[-\frac{(E-E_0)^2}{2\sigma_E^2}\Big] 
\prod_{k=2}^{d} 
\frac{1}{\sqrt{2\pi}\, \sigma_k} \exp\Big[-\frac{(x_k-x_{k,0})^2}{2\sigma_k^2}\Big]
\end{equation}
of a monoenergetic proton beam with energy $E_0$. The isocenter of the beam
consisting of $\psi_0$ protons
is located at 
the point $\x_0 = (0, x_{2,0},\dots,  x_{d,0})^\top\in\Gamma_{\mathrm{in}}$.
We choose the energy spread $ \sigma_E =0.01\, E_0$ and the spatial standard deviation $\sigma_k = 0.3$ for all $k\geq 2$~\cite{Stammer2025}.
In one space dimension, \eqref{eq:bdrcon0} reduces to a Gaussian distribution in energy.
The corresponding first moment is prescribed as in~\cite{pichard2016}
\begin{equation}\label{eq:bdrcon1}
    \hat{\boldsymbol{\psi}}^{(1)}(\x,E) = f\,\hat \psi^{(0)}(\x, E)\,\mathbf{e}_1, \qquad f = 0.9999,
\end{equation}
where $\mathbf{e}_1$ is the unit vector in the $x_1$-direction. This choice
approximates the moments of a nearly perfectly collimated particle beam.
We set $E_{\max} = 1.1\, E_0$ to ensure that the full energy spread of the beam is captured.

The initial fluence $\psi$ at $E_{\max}$ is usually assumed to be zero~\cite{frank2007}.
However, the angular moments associated with the trivial angular distribution lie on the boundary of the realizable set. 
To avoid this, we prescribe the initial condition~\cite{pichard2016}
\begin{equation}\label{eq:initcon}
    \pz(E_{\max}) \equiv 10^{-15}, \qquad \po(E_{\max}) \equiv0
\end{equation}
in the whole domain $\D$. 
Multiplying~\eqref{eq:initcon} by $S(E_{\max})$, we obtain the initial condition for $(Su)$.

During the energy stepping, the dose is accumulated using the composite trapezoidal rule, i.e.,
\begin{equation*}
    \int_{E_{\min}}^{E_{\max}} (S\pz)_i\,\diff E  \approx \sum_{n = 1}^{N-1} \frac{(S\pz)_i^n + (S\pz)_i^{n+1}}{2}(E^{n+1}-E^n).
\end{equation*}
The integral over the remaining energy interval $[0,E_{\min}]$ is treated differently.
We assume that the proton fluence at energy $E_{\min}$ deposits its residual energy locally.
Under this assumption, the stopping power in the interval $[0,E_{\min}]$ is replaced by
\begin{equation}\label{eq:S0}
    S_0 = \frac{E_{\min}}{R(E_{\min})},
\end{equation}
where $R(E)$ is the proton range given by \eqref{eq:range}.
This treatment is similar to the fictitious group assumption in the context of a multigroup method; see, e.g.,~\cite{morel1981}.
The choice~\eqref{eq:S0} is consistent with the interpretation of the stopping power as the energy loss per unit path length. Thus, we approximate the dose by 
\begin{equation*}
    D(\x) \approx D_h(\x) = \sum_{i= 1}^{N_h} D_i\varphi_i(\x),
\end{equation*}
where
\begin{equation*}
    D(\x_i)\approx D_i = \frac{1}{\rho(\x_i)}\left(S_0 \psi^{(0),1}_i E_{\min} + \sum_{n = 1}^{N-1} \frac{(S\pz)_i^n + (S\pz)_i^{n+1}}{2}(E^{n+1}-E^n)\right)
\end{equation*}
and $\psi^{(0),1}_i = (S\pz)_i^{1}/S_i(E_{\min})$.

\section{Numerical examples}
\label{sec:examples}

To assess the proposed limiting strategy in the context of dose calculation, we apply our realizability-preserving MCL scheme to representative test problems.
For the energy discretization, we use Heun's method, a second-order explicit SSP-RK scheme.
In view of condition \eqref{eq:CFL}, the time step $\Delta E$ is determined using the formula~\cite{guermond2016, kuzmin2020, kuzmin2023} 
\begin{equation*}
    \Delta E = \frac{\mathrm{CFL}}{\max_{i\in\{1,\ldots,N_h\}}\frac{2}{m_i S_i }\sum_{j\in\N^*_i}d_{ij}},
\end{equation*}
where $\mathrm{CFL} \leq 1$ is a given threshold.
This choice of $\Delta E$ guarantees realizability, as shown in Sections~\ref{sec:loworder} and~\ref{sec:MCL}.

The implementation of MCL used in our numerical experiments is based on the open-source \texttt{C++} finite element library MFEM~\cite{anderson2021,andrej2024,mfem}. 
The two- and three-dimensional results are visualized in Paraview~\cite{ayachit2015}.

\subsection{Analytical model}

We begin by validating the proposed method against an analytical reference solution derived and verified with established Monte Carlo dose engines in~\cite{ashby2025}. Neglecting the scattering effects and assuming a perfectly collimated beam propagation in a single direction through a homogeneous medium reduces the Fokker--Planck equation~\eqref{eq:fp} to a linear one-dimensional transport equation, which can be solved via the method of characteristics.
For a prescribed inflow fluence $g(E)$, the corresponding dose is obtained by integrating the analytical solution weighted by the stopping power~\eqref{eq:SBK}, yielding~\cite{ashby2025}
\begin{equation}\label{eq:refSol}
    \begin{split}
        D_\mathrm{ref}(x) = \frac{1}{\beta p \rho(x)} \int_0^{E_{\max}}\left(E^p + \frac{x}{\beta}\right)^{\frac{1-p}{p}}g\left(\left(E^p + \frac{x}{\beta} \right)^{\frac{1}{p}} \right)\,\diff E.
    \end{split}
\end{equation}
In this first test, we set the scattering power to $T = 0$ for a direct comparison with the analytical reference solution.
Note that the scattering step~\eqref{eq:scatteringstep} reduces to multiplication by unity.
We consider a monoenergetic proton beam with energy $E_0= 62\mev$ consisting of $\psi_0 = 1.21\times 10^9$ protons in a water phantom.
Thus, the material properties are $\rho= 1$, $\beta = 0.0022$, and $p = 1.77$ throughout the computational domain $\D = [0,4]\cm$.
The beam is prescribed at $x= 0$ using~\eqref{eq:bdrcon0}--\eqref{eq:bdrcon1}.
Consistently, the inflow function used in the reference solution~\eqref{eq:refSol} is given by
\begin{equation*}
    g(E) = \psi_0 \,\frac{1}{\sqrt{2\pi}\, \sigma_E} \exp\Big[-\frac{(E-E_0)^2}{2\sigma_E^2}\Big].
\end{equation*}
In two further tests, we perform simulations in the three-dimensional domain $\D = [0,4]\cm \times [0, 1.5]\cm \times [0, 1.5]\cm $. To assess the accuracy of the dose distributions computed with and without scattering, we compare them with the reference solution~\eqref{eq:refSol}. For this comparison, the dose is integrated over the $y$-$z$ plane perpendicular to the direction of beam propagation, as done in~\cite{ashby2025}.

\begin{figure}[h!]
    \centering
    \begin{subfigure}[b]{0.48\textwidth}
        \centering
        \includegraphics[width=\textwidth, trim=1cm 6cm 2cm 7cm, clip]{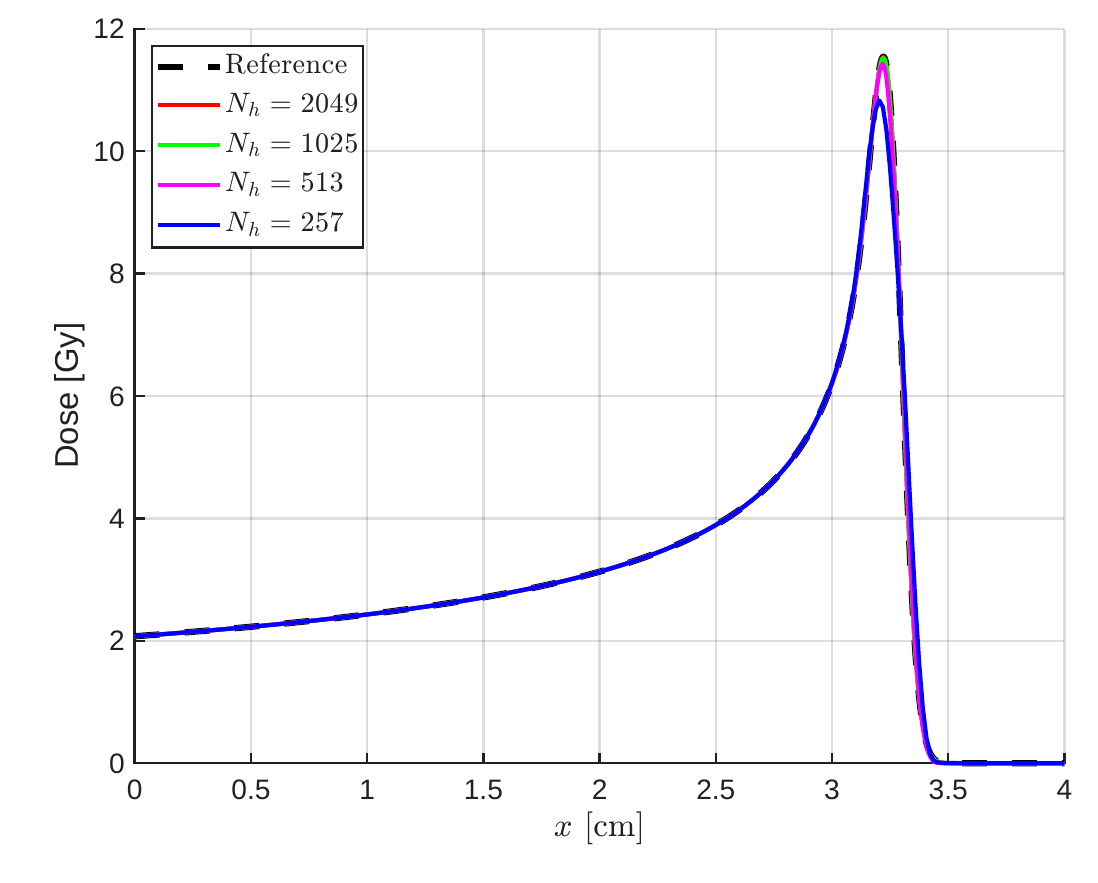}
        \caption{One-dimensional dose distributions obtained on meshes of varying resolution.}
        \label{fig:1ddose}
    \end{subfigure}
    \hfill
    \begin{subfigure}[b]{0.48\textwidth}
        \centering
        \includegraphics[width=\textwidth, trim=1cm 6cm 2cm 7cm, clip]{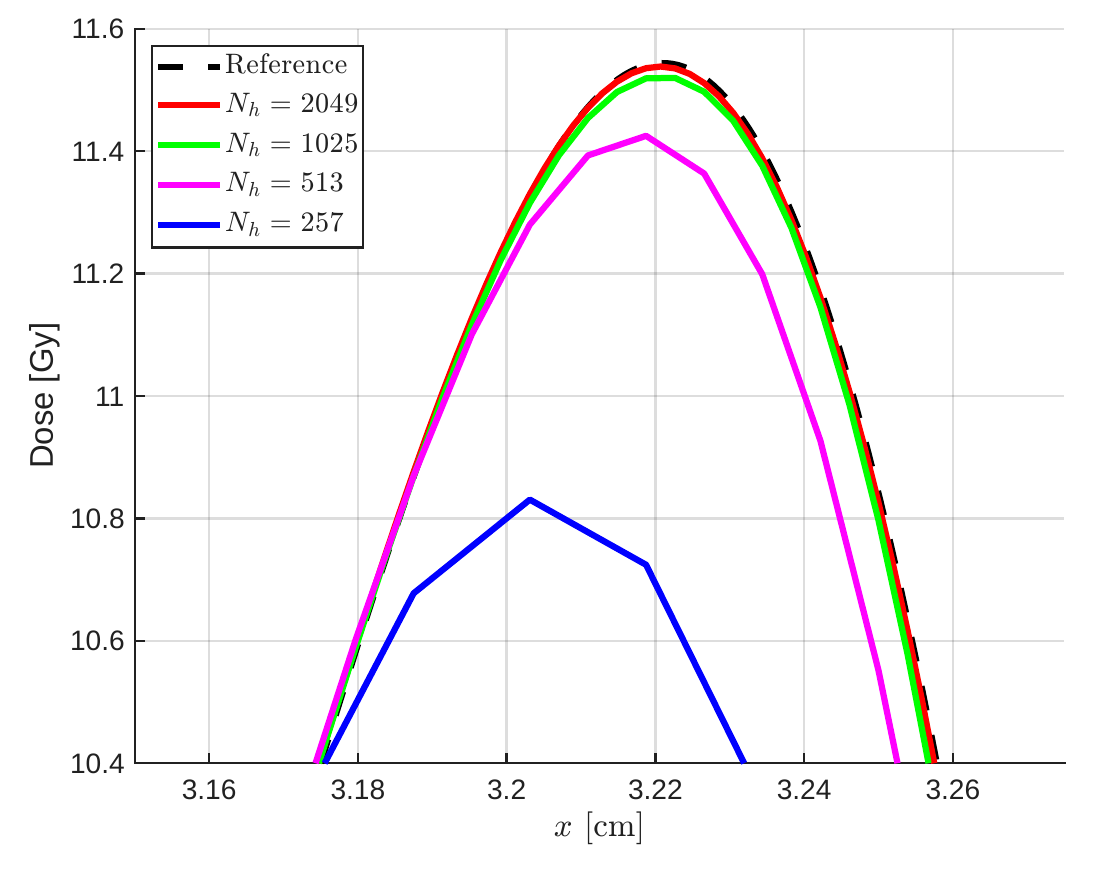}
        \caption{Close up of Bragg-peak region for the one-dimensional dose distributions.}
        \label{fig:1dBraggpeak}
    \end{subfigure}
    \caption{One-dimensional dose distribution of a $62\mev$ proton beam in a water phantom computed with the proposed MCL scheme on meshes with $N_h\in\{257, 513, 1025, 2049\}$ and $\mathrm{CFL}=0.5$. Scattering is neglected.
    Results are compared with the reference solution~\eqref{eq:refSol}.}
    \label{fig:refsolcomparison}
\end{figure}
Figure~\ref{fig:refsolcomparison} shows the one-dimensional dose distributions computed with the proposed MCL scheme for the $M_1$ model on a hierarchy of meshes.
Away from the Bragg peak, even coarse meshes provide accurate approximations of the dose profile, when compared to the reference solution~\eqref{eq:refSol}.
However, near the Bragg peak, coarse discretizations exhibit noticeable peak clipping. 
As the mesh is refined, the resolution of the peak improves significantly.
On the finest mesh with $N_h=2049$ nodes, the numerical solution agrees almost perfectly with the reference solution, despite employing the reduced $M_1$ moment model.

\begin{figure}[h!]
    \centering
    \begin{subfigure}[b]{\textwidth}
        \centering
        \includegraphics[width=\textwidth, trim=0cm 2.5cm 0cm 0cm, clip]{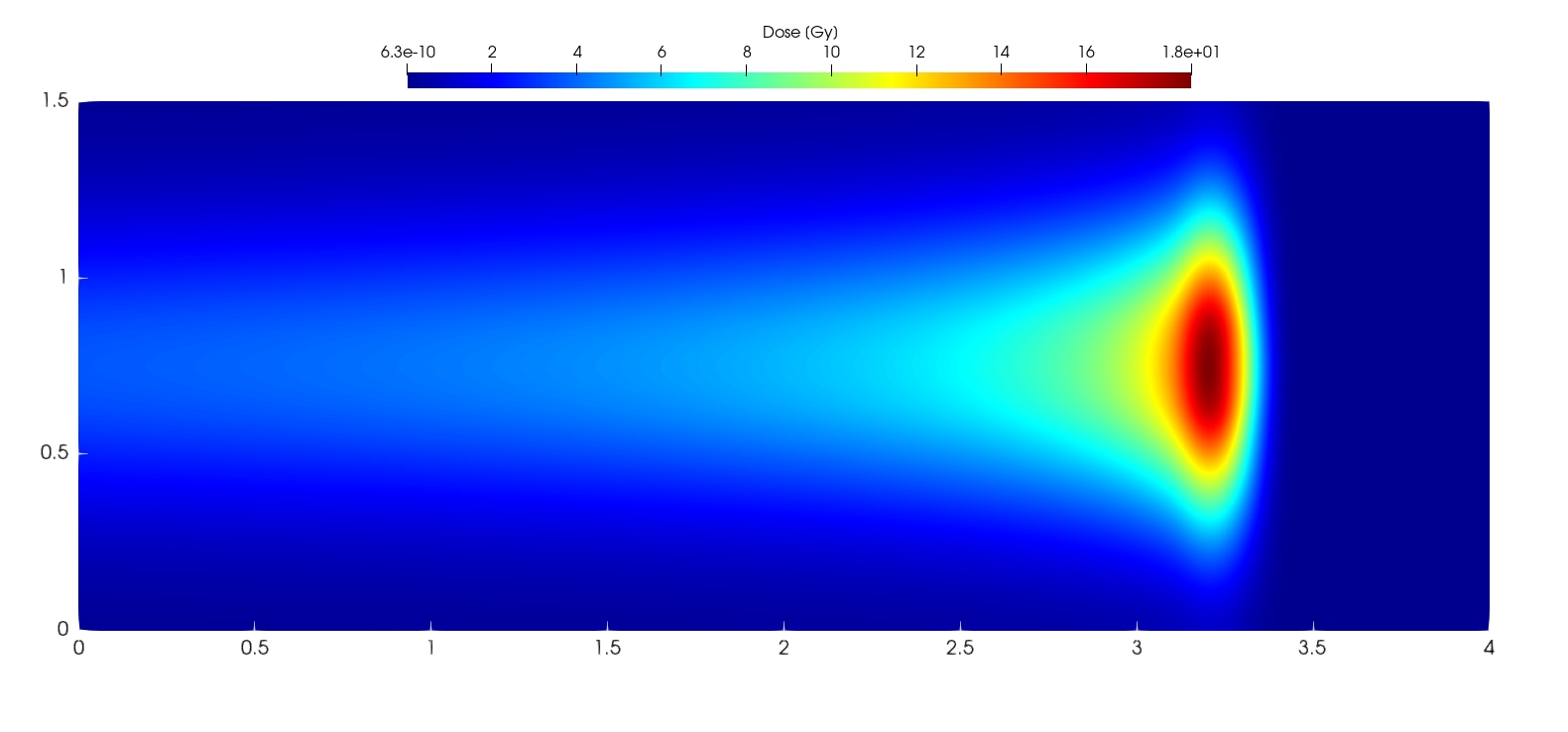}
        \label{fig:3dclosedsol}
    \end{subfigure}
    \hfill
    \begin{subfigure}[b]{\textwidth}
        \centering
        \includegraphics[width=\textwidth, trim=0cm 2.5cm 0cm 3.2cm, clip]{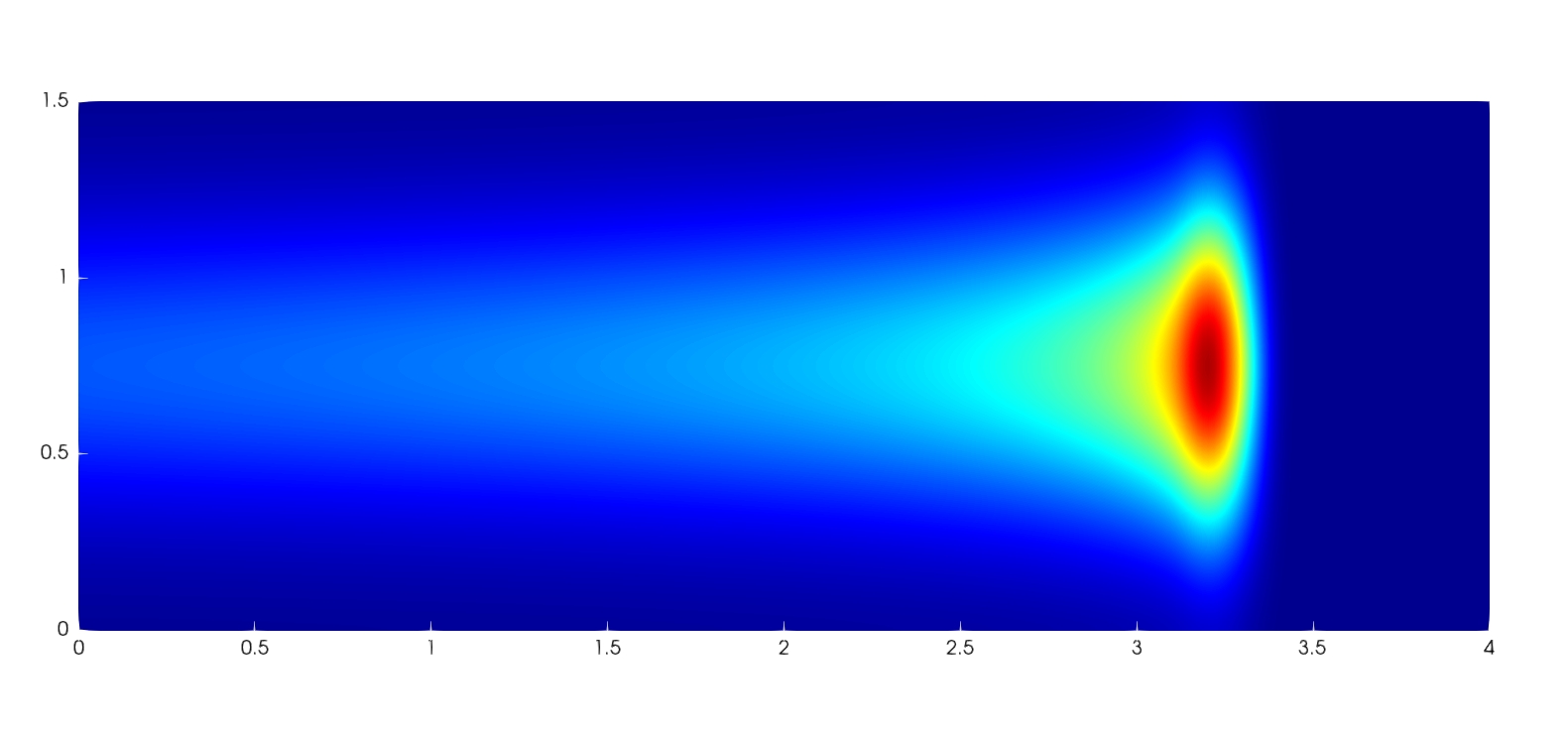}
        \label{fig:3dclosedsolT>0}
    \end{subfigure}
    \caption{Three-dimensional dose distribution of a $62\mev$ proton beam in a water phantom computed with the proposed realizability-preserving MCL scheme on a uniform hexahedral mesh with $N_h = 257 \times 97 \times 97$ nodes and $\mathrm{CFL}=0.5$. Shown are slices at $z = 0.75\cm$ without scattering (top) and with scattering (bottom).}
    \label{fig:3D_closedform_slice}
\end{figure}

\begin{figure}[h!]
    \centering
    \begin{subfigure}[b]{0.48\textwidth}
        \centering
        \includegraphics[width=\textwidth, trim=1cm 6cm 2cm 7cm, clip]{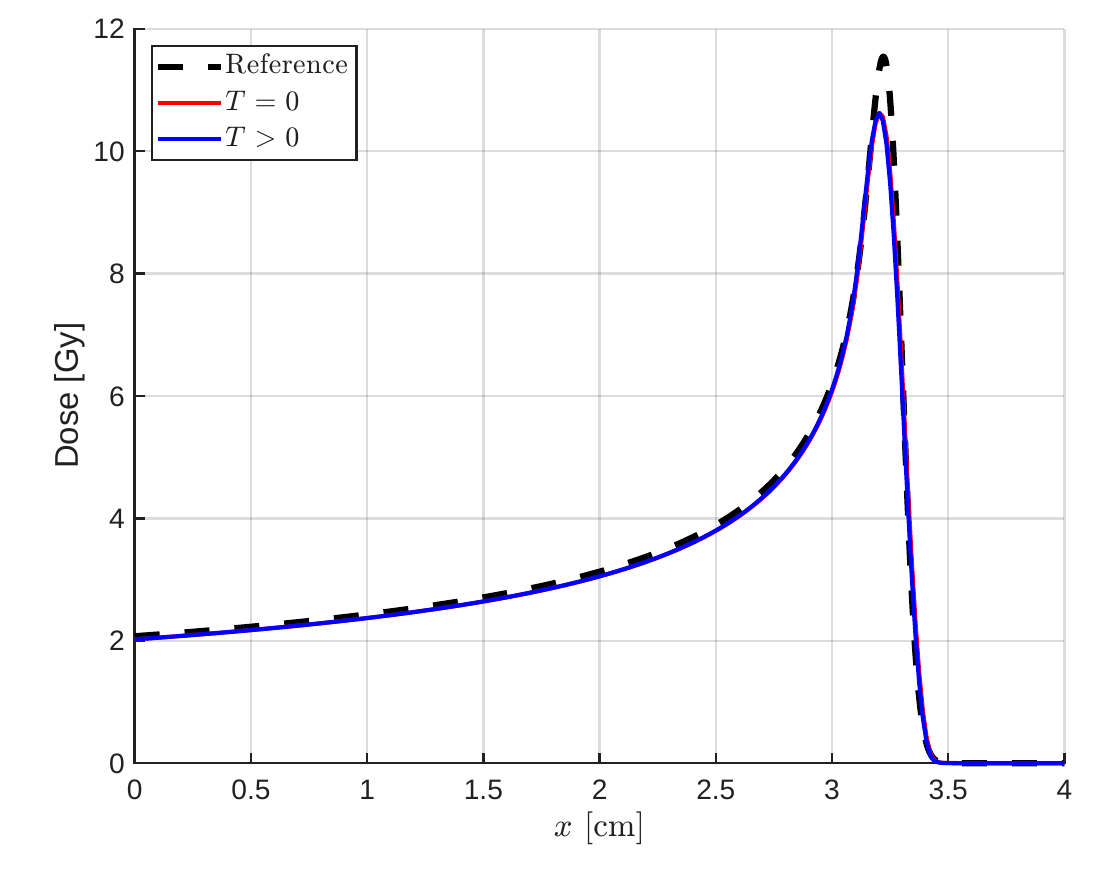}
        \caption{Three-dimensional dose distributions integrated over $y$-$z$ planes.}
        \label{fig:3ddose}
    \end{subfigure}
    \hfill
    \begin{subfigure}[b]{0.48\textwidth}
        \centering
        \includegraphics[width=\textwidth, trim=1cm 6cm 2cm 7cm, clip]{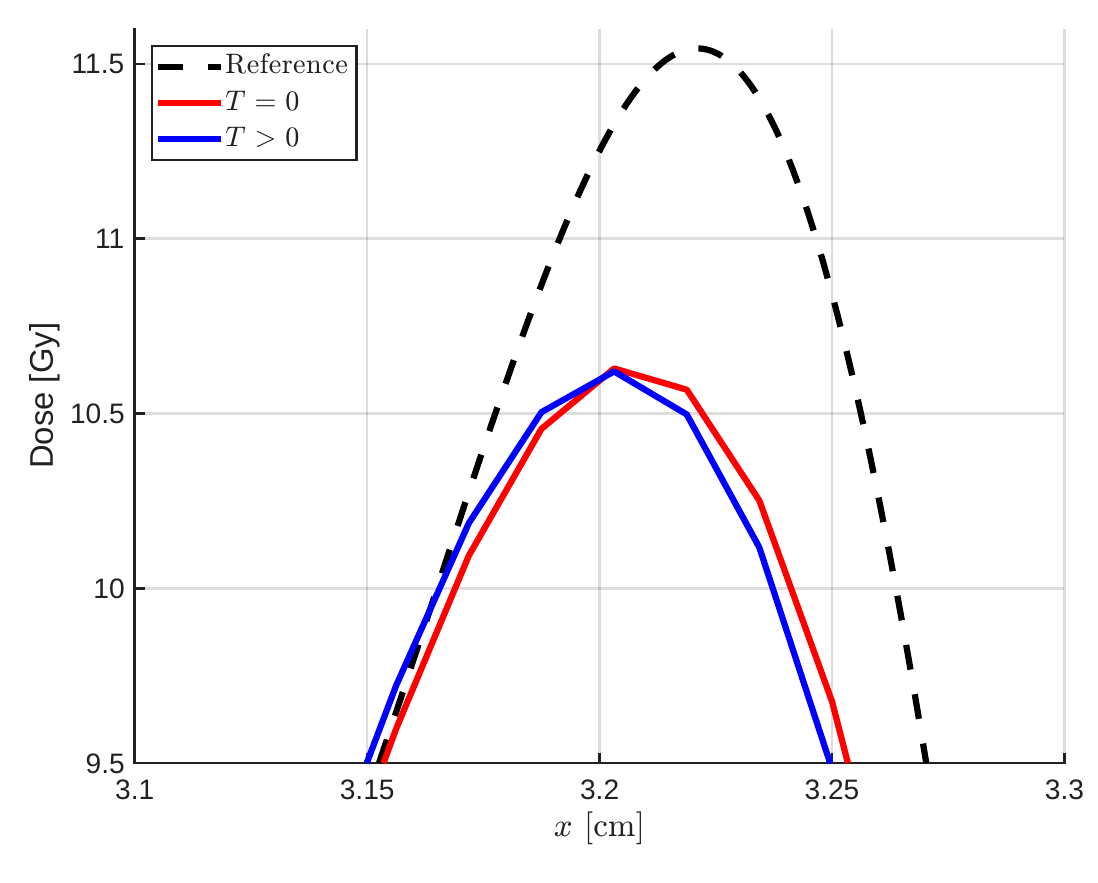}
        \caption{Close up of Bragg-peak region for the three-dimensional dose distributions integrated over $y$-$z$ planes.}
        \label{fig:3dBraggpeak}
    \end{subfigure}
    \caption{Three-dimensional dose distributions of a $62\mev$ proton beam in a water phantom computed with the proposed MCL scheme on a mesh with $N_h = 257\times 97 \times 97$ nodes and $\mathrm{CFL}=0.5$ with and without scattering effects.
    Results are integrated over $y$-$z$ planes and compared with the reference solution~\eqref{eq:refSol}.}
    \label{fig:3Dintegrated}
\end{figure}

Slices of the three-dimensional solutions with and without scattering at $z = 0.75\cm$ are shown in Figure~\ref{fig:3D_closedform_slice}.
Both simulations yield physically consistent dose distributions without any visible numerical artifacts. 
In the Bragg-peak region, the solution including scattering exhibits a slightly reduced peak.

The corresponding depth-dose curves obtained by integration over the $y$-$z$ plane are shown in Figure~\ref{fig:3Dintegrated}.
Both curves compare well with the reference solution.
Minor peak clipping is observed, which can be attributed to the spatial resolution in the $x$-direction.
The nearly identical peak magnitudes in the integrated curves indicate that the slight difference in peak dose observed in Figure~\ref{fig:3D_closedform_slice} is due to lateral spreading of the dose profile caused by scattering.

\subsection{65 MeV beam in a patient}
To illustrate the performance of the proposed limiting techniques in the presence of material discontinuities, we consider a three-dimensional simulation of a $65\mev$ proton beam consisting of $\psi_0 = 1.21\times 10^9$ protons propagating in the $x$-direction, prescribed at $x = 0$ using~\eqref{eq:bdrcon0}--\eqref{eq:bdrcon1}.
The computational domain is $\D = [0,4]\cm \times [0, 1.5]\cm \times [0, 1.5]\cm$, which is decomposed into four slabs of muscle, bone, lung, and water,
\begin{equation*}
\begin{split}
\D_{\mathrm{muscle}}&=\{(x,y,z)\in\D\,:\, 0 \le x < 1 \},\\
\D_{\mathrm{bone}}&=\{ (x,y,z) \in \D \,:\, 1 \le x < 1.25 \},\\
\D_{\mathrm{lung}}&=\{ (x,y,z) \in \D \,:\, 1.25 \le x < 3 \},\\
\D_{\mathrm{water}}&=\{(x,y,z)\in\D\,:\,3\leq x\leq 4\},
\end{split}
\end{equation*}
respectively.
The stopping and scattering powers are set constant in space within each material slab using the parameters given in Table~\ref{tab:BK_XS}.

\begin{figure}[h!]
    \centering
    \begin{subfigure}[b]{\textwidth}
        \centering
        \includegraphics[width=\textwidth, trim=0cm 2.5cm 0cm 0cm, clip]{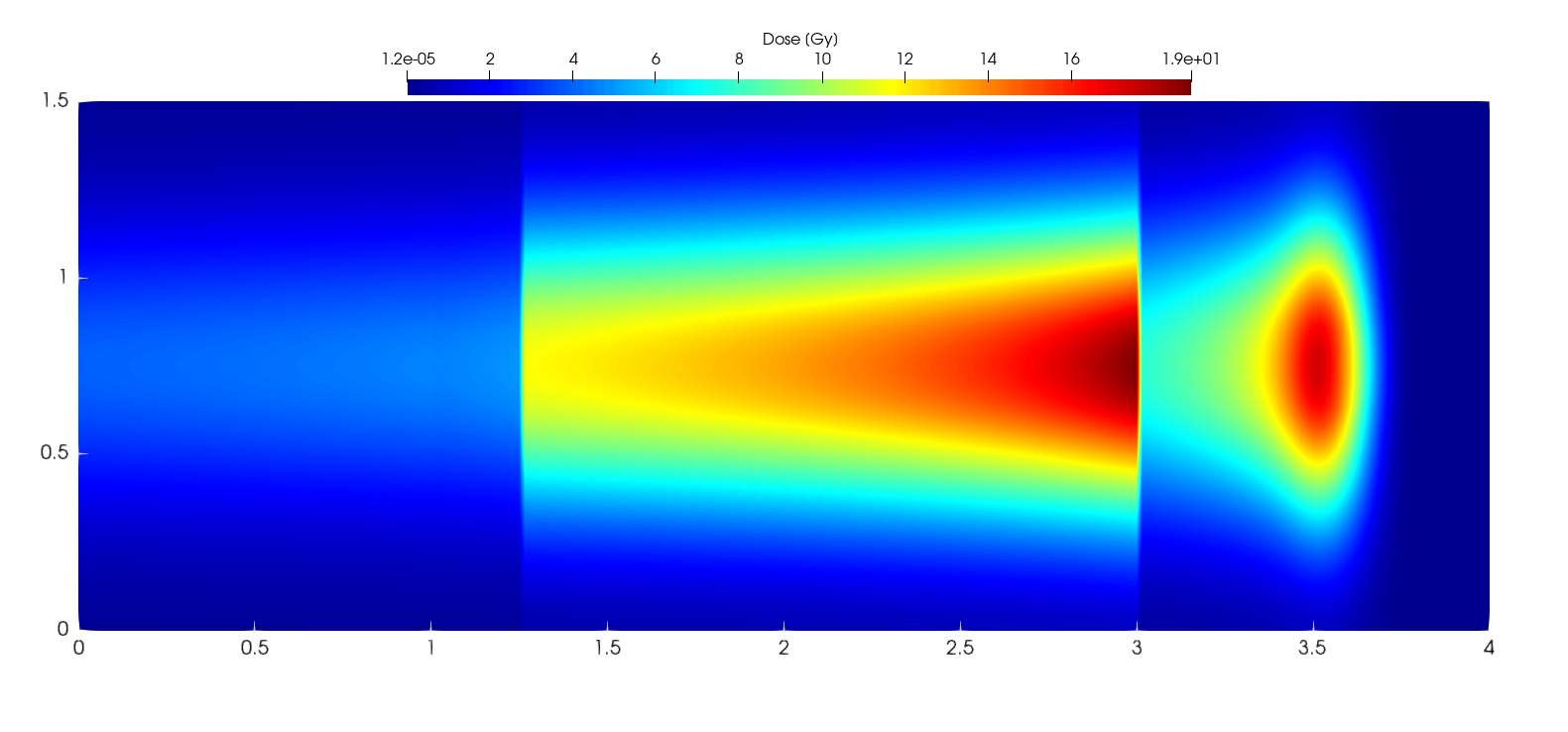}
        \caption{Dose.}
        \label{fig:3DmultimatDose}
    \end{subfigure}
    \hfill
    \begin{subfigure}[b]{\textwidth}
        \centering
        \includegraphics[width=\textwidth, trim=0cm 2.5cm 0cm 0cm, clip]{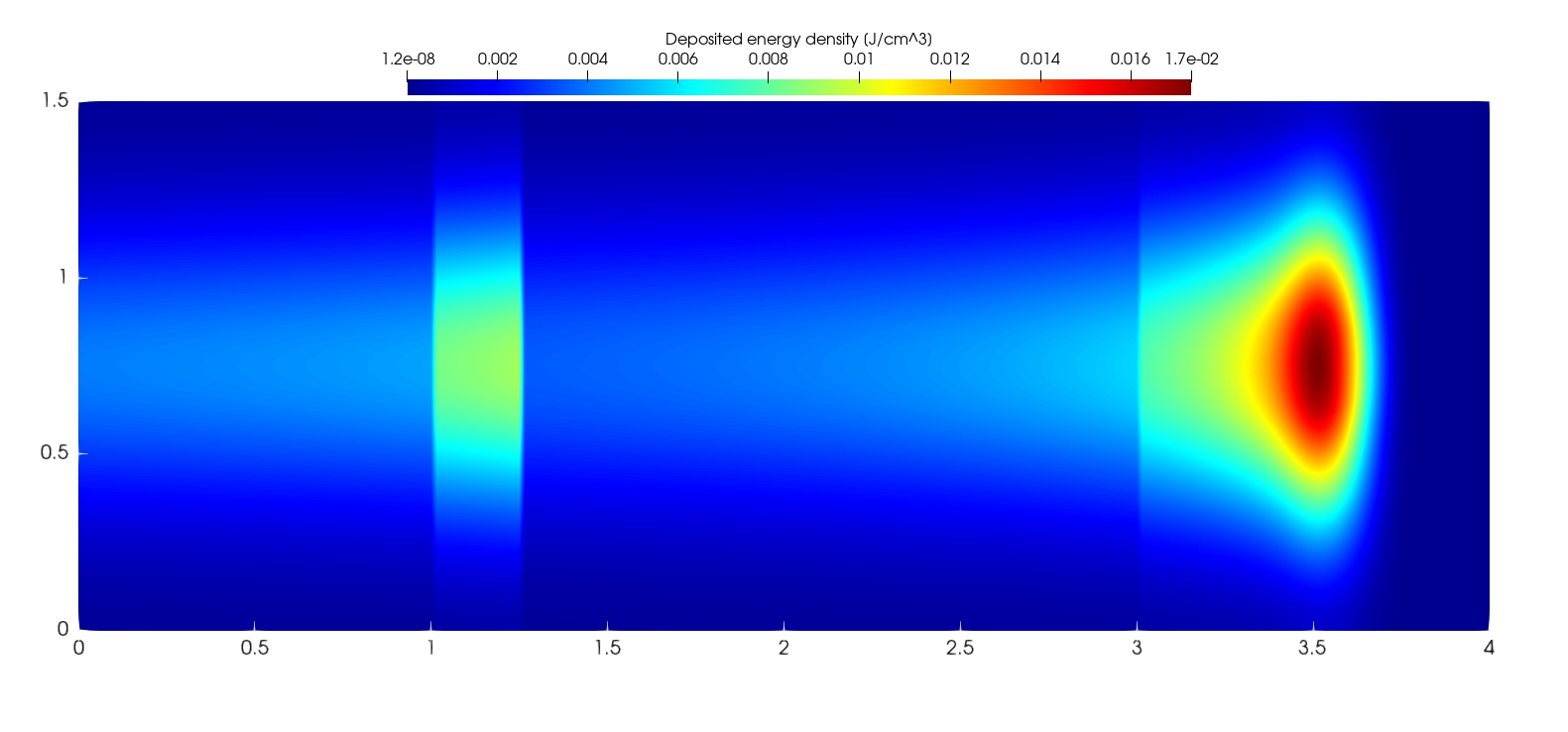}
        \caption{Deposited energy density.}
        \label{fig:3DmultimatEnergy}
    \end{subfigure}
    \caption{Three-dimensional dose distribution of a $65\mev$ proton beam in a heterogeneous multi-material geometry computed with the proposed MCL scheme on a uniform hexahedral mesh with $N_h = 257 \times 97 \times 97$ nodes and $\mathrm{CFL}=0.5$. Shown are slices at $z = 0.75\cm$ with dose (top) and deposited energy density (bottom).}
    \label{fig:3Dmultimat}    
\end{figure}

Figure~\ref{fig:3DmultimatDose} shows a slice of the three-dimensional dose distribution at $z = 0.75\cm$.
Since the interface between muscle and bone is not clearly visible in the dose distribution, we additionally visualize the deposited energy density, $\rho D$, in Figure~\ref{fig:3DmultimatEnergy}.
All material interfaces are well resolved without any oscillations.
No nonphysical states were detected during the simulation and no numerical artifacts are visible.
This result demonstrates that the proposed realizability-preserving MCL scheme can resolve material discontinuities sharply while maintaining physical consistency. 

\subsection{Double beam problem}
\label{sec:doublebeam}
In the final numerical example, we illustrate a well-known drawback of the $M_1$ model.
It turns out that the $M_1$ model cannot distinguish between two overlapping beams.
This limitation arises because the zeroth and first angular moments of two intersecting beams coincide with those of a single beam propagating in the mean direction~\cite[Examples~2 and~3]{pichard2017}.
To illustrate this effect, we consider the two-dimensional computational domain $\D = [0,4]\cm\times[0,4]\cm$ consisting of water.
Two proton beams with energy $E_0 = 62\mev$, each consisting of $\psi_0 = 1.21\times 10^9$ protons, are prescribed at the midpoints of the left and lower boundaries, propagating in $x$- and $y$-directions, respectively.

\begin{figure}[h!]
	\centering
		\includegraphics[trim={14cm 0cm 10cm 0cm},clip, width = 0.6\textwidth]{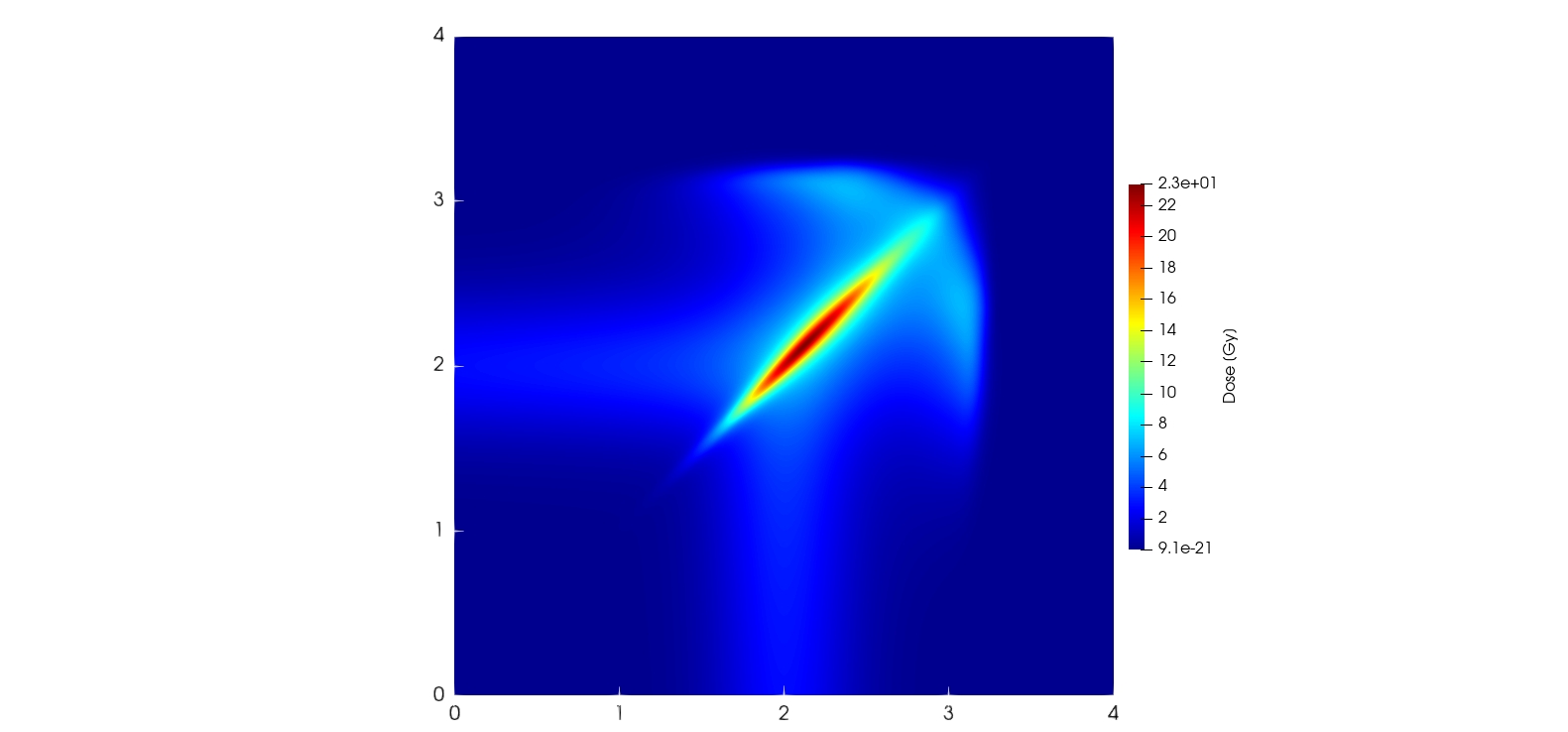}
		\caption{Dose distribution for two perpendicular proton beams in a water phantom computed with the proposed MCL scheme on a uniform rectangular mesh with $N_h = 257\times 257$ nodes and $\mathrm{CFL}= 0.5$.}
        \label{fig:2beams}
\end{figure}

As shown in Fig.~\ref{fig:2beams}, the beams intersect near the center of the computational domain and merge into a single beam traveling along the diagonal.
This behavior is consistent with the literature and with \cite[Examples 2 and 3]{pichard2017}. 
Again, no nonphysical states were detected throughout the simulation and no numerical instabilities are visible in the numerical solution.


\section{Conclusions}
\label{sec:conclusion}

We have proposed a realizability-preserving MCL scheme for continuous finite element discretizations of the energy-dependent $M_1$ moment model of proton transport for dose calculation. Energy is treated as a pseudo-time variable, and a Strang splitting approach is employed to handle the scattering-induced forcing terms. The transport subproblem is discretized in energy using an explicit SSP-RK method. In this way, the moments are evolved backward in energy while preserving realizability at each substep. The dose is computed during the energy evolution using the trapezoidal rule to approximate the integral of the zeroth moment weighted by the stopping power.

Numerical experiments confirm that the proposed scheme produces stable and physically consistent dose distributions for both homogeneous media and heterogeneous material slabs in single-beam scenarios. The Bragg peak is well approximated when the mesh size is sufficiently small, and the scheme captures material interfaces in a sharp and well-resolved manner. However, the double-beam problem in Section~\ref{sec:doublebeam} highlights a fundamental limitation of the $M_1$ moment model: overlapping beams merge into a single beam propagating in the mean direction.

This observation motivates the extension of the MCL methodology to continuous Galerkin discretizations of the $M_2$ moment model, in which the second moment $\pt$ is computed explicitly and the third moment is modeled by a closure relation~\cite{pichard2017}. In this setting, additional admissibility conditions arise in the form of eigenvalue constraints on the second-moment tensor. Limiting frameworks for tensor fields have been proposed in~\cite{kuzmin2020,lohmann2017b,lohmann2019}.

\bibliographystyle{plain}
\bibliography{bibliography}

@Article{anderson2021,
  author  = {Anderson, Robert and Andrej, Julian and Barker, Andrew and Bramwell, Jamie and Camier, Jean-Sylvain and Cerveny, Jakub and Dobrev, Veselin and Dudouit, Yohann and Fisher, Aaron and Kolev, Tzanio and Pazner, Will and Stowell, Mark and Tomov, Vladimir and Akkerman, Ido and Dahm, Johann and Medina, David and Zampini, Stefano},
  journal = {Comput. Math. Appl.},
  title   = {{MFEM: A} modular finite element methods library},
  year    = {2021},
  pages   = {42--74},
  volume  = {81},
  doi     = {10.1016/j.camwa.2020.06.009},
}

@Book{ayachit2015,
  author    = {Ayachit, Utkarsh},
  publisher = {Kitware, Inc.},
  title     = {{The ParaView guide: A Parallel Visualization Application}},
  year      = {2015},
  isbn      = {9781930934306},
}

@Article{barrenechea2017b,
  author  = {Barrenechea, Gabriel R. and Knobloch, Petr},
  journal = {Appl. Numer. Math.},
  title   = {Analysis of a group finite element formulation},
  year    = {2017},
  pages   = {238--248},
  volume  = {118},
  doi     = {10.1016/j.apnum.2017.03.008},
}

@article{burchard2003,
  title={A high-order conservative {P}atankar-type discretisation for stiff systems of production--destruction equations},
  author={Burchard, Hans and Deleersnijder, Eric and Meister, Andreas},
  journal={Appl. Numer. Math.},
  volume={47},
  number={1},
  pages={1--30},
  year={2003},
  doi={10.1016/S0168-9274(03)00101-6}
}

@article{fletcher1983,
title = {The group finite element formulation},
journal = {Computer Methods Appl. Mech. Engrg.},
volume = {37},
number = {2},
pages = {225--244},
year = {1983},
doi = {10.1016/0045-7825(83)90122-6},
author = {C.A.J. Fletcher}
}

@Article{guermond2016,
  author    = {Guermond, Jean-Luc and Popov, Bojan},
  journal   = {SIAM J. Numer. Anal.},
  title     = {Invariant domains and first-order continuous finite element approximation for hyperbolic systems},
  year      = {2016},
  issn      = {1095-7170},
  month     = jan,
  number    = {4},
  pages     = {2466--2489},
  volume    = {54},
  doi       = {10.1137/16M1074291},
  publisher = {Society for Industrial & Applied Mathematics (SIAM)},
}

@Article{hajduk2021,
  author    = {Hennes Hajduk},
  journal   = {Comput. Math. Appl.},
  title     = {Monolithic convex limiting in discontinuous {G}alerkin discretizations of hyperbolic conservation laws},
  year      = {2021},
  issn      = {0898-1221},
  month     = apr,
  pages     = {120--138},
  volume    = {87},
  doi       = {10.1016/j.camwa.2021.02.012},
  publisher = {Elsevier BV},
}

@Article{kuzmin2010a,
  author    = {Dmitri Kuzmin and Matthias M\"oller and John N. Shadid and Mikhail Shashkov},
  journal   = {J. Comput. Phys.},
  title     = {Failsafe flux limiting and constrained data projections for equations of gas dynamics},
  year      = {2010},
  issn      = {0021-9991},
  month     = nov,
  number    = {23},
  pages     = {8766--8779},
  volume    = {229},
  doi       = {10.1016/j.jcp.2010.08.009},
  publisher = {Elsevier BV},
}

@Article{kuzmin2020,
  author    = {Dmitri Kuzmin},
  journal   = {Computer Methods Appl. Mech. Engrg.},
  title     = {Monolithic convex limiting for continuous finite element discretizations of hyperbolic conservation laws},
  year      = {2020},
  issn      = {0045-7825},
  month     = apr,
  pages     = {112804},
  volume    = {361},
  doi       = {10.1016/j.cma.2019.112804},
  publisher = {Elsevier BV},
}

@Book{kuzmin2023,
  author    = {Dmitri Kuzmin and Hennes Hajduk},
  publisher = {World Scientific},
  title     = {{Property-Preserving Numerical Schemes for Conservation Laws}},
  year      = {2023},
  doi       = {10.1142/13466},
}

@Article{lohmann2017b,
  author    = {Christoph Lohmann},
  journal   = {J. Comput. Phys.},
  title     = {Flux-corrected transport algorithms preserving the eigenvalue range of symmetric tensor quantities},
  year      = {2017},
  issn      = {0021-9991},
  month     = dec,
  pages     = {907-926},
  volume    = {350},
  doi       = {10.1016/j.jcp.2017.09.009},
  publisher = {Elsevier BV},
}

@Book{lohmann2019,
  author    = {Lohmann, Christoph},
  publisher = {Springer Spektrum},
  title     = {{Physics-Compatible Finite Element Methods for Scalar and Tensorial Advection Problems}},
  year      = {2019},
  isbn      = {9783658277376},
  doi       = {10.1007/978-3-658-27737-6},
}

@misc{mfem,
key = {MFEM},
title = {{MFEM:} {M}odular {F}inite {E}lement {M}ethods [{S}oftware]},
howpublished = {\url{https://mfem.org}},
}

@Book{kuzmin2014a,
  author    = {Kuzmin, Dmitri and H{\"a}m{\"a}l{\"a}inen, Jari},
  publisher = {SIAM},
  title     = {{Finite Element Methods for Computational Fluid Dynamics: A Practical Guide}},
  year      = {2014},
}

@Article{alldredge2012,
  author    = {Alldredge, Graham W. and Hauck, Cory D. and Tits, Andre L.},
  journal   = {SIAM J. Sci. Comput.},
  title     = {High-order entropy-based closures for linear transport in slab geometry {II}: {A} computational study of the optimization problem},
  year      = {2012},
  number    = {4},
  pages     = {B361--B391},
  volume    = {34},
  doi       = {10.1137/11084772X},
  publisher = {SIAM},
}

@Article{andrej2024,
  author  = {Andrej, Julian and Atallah, Nabil and B{\"a}cker, Jan-Phillip and Camier, Jean-Sylvain and Copeland, Dylan and Dobrev, Veselin and Dudouit, Yohann and Duswald, Tobias and Keith, Brendan and Kim, Dohyun and others},
  journal = {Int. J. High Perform. Comput. Appl.},
  title   = {High-performance finite elements with {MFEM}},
  year    = {2024},
  number  = {5},
  pages   = {447--467},
  volume  = {38},
  doi     = {10.1177/10943420241261981},
}

@Article{berthon2007,
  author    = {Berthon, Christophe and Charrier, Pierre and Dubroca, Bruno},
  journal   = {J. Sci. Comput.},
  title     = {An {HLLC} scheme to solve the {M}1 model of radiative transfer in two space dimensions},
  year      = {2007},
  pages     = {347--389},
  volume    = {31},
  publisher = {Springer},
  url       = {https://link.springer.com/article/10.1007/s10915-006-9108-6},
}

@Article{brunner2000,
  author  = {Brunner, Thomas A. and Holloway, James P.},
  journal = {Transport},
  title   = {Two new boundary conditions for use with the maximum entropy closure and an approximate {R}iemann solver},
  year    = {2000},
  pages   = {3},
  volume  = {10},
  url     = {https://www.researchgate.net/profile/James-Holloway-3/publication/268419497_TWO_NEW_BOUNDARY_CONDITIONS_FOR_USE_WITH_THE_MAXIMUM_ENTROPY_CLOSURE_AND_AN_APPROXIMATE_RIEMANN_SOLVER/links/552fb04d0cf2acd38cbc5174/TWO-NEW-BOUNDARY-CONDITIONS-FOR-USE-WITH-THE-MAXIMUM-ENTROPY-CLOSURE-AND-AN-APPROXIMATE-RIEMANN-SOLVER.pdf},
}

@Article{brunner2001,
  author    = {Brunner, Thomas A. and Holloway, James Paul},
  journal   = {J. Quant. Spectrosc. Radiat. Transfer},
  title     = {One-dimensional {R}iemann solvers and the maximum entropy closure},
  year      = {2001},
  number    = {5},
  pages     = {543--566},
  volume    = {69},
  doi       = {10.1016/s0022-4073(00)00099-6},
  publisher = {Elsevier},
}

@Article{chidyagwai2018,
  author    = {Chidyagwai, Prince and Frank, Martin and Schneider, Florian and Seibold, Benjamin},
  journal   = {J. Comput. Appl. Math.},
  title     = {A comparative study of limiting strategies in discontinuous {G}alerkin schemes for the {$M_1$} model of radiation transport},
  year      = {2018},
  pages     = {399--418},
  volume    = {342},
  doi       = {10.1016/j.cam.2018.04.017},
  publisher = {Elsevier},
}

@Article{coulombel2006,
  author    = {Coulombel, Jean-Fran{\c{c}}ois and Goudon, Thierry},
  journal   = {J. Hyperbol. Differ. Eq.},
  title     = {Entropy-based moment closure for kinetic equations: {R}iemann problem and invariant regions},
  year      = {2006},
  issn      = {1793-6993},
  month     = dec,
  number    = {04},
  pages     = {649--671},
  volume    = {3},
  doi       = {10.1142/S0219891606000951},
  publisher = {World Scientific},
}

@Article{frank2012,
  author  = {Frank, Martin and Hauck, Cory D. and Olbrant, Edgar},
  journal = {arXiv preprint arXiv:1208.0772},
  title   = {Perturbed, entropy-based closure for radiative transfer},
  year    = {2012},
  doi     = {10.48550/arXiv.1208.0772},
}

@Article{hauck2011,
  author  = {Hauck, Cory D.},
  journal = {Commun. Math. Sci},
  title   = {High-order entropy-based closures for linear transport in slab geometry},
  year    = {2011},
  number  = {1},
  pages   = {187--205},
  volume  = {9},
  url     = {https://www.math.umd.edu/~tadmor/ki_net/pubs/files/FRG-2010-Hauck-Cory.entropy_kinetic.pdf},
}

@TechReport{kershaw1976,
  author      = {Kershaw, David S.},
  institution = {Lawrence Livermore National Lab.(LLNL), Livermore, CA (United States)},
  title       = {Flux limiting nature`s own way -- {A} new method for numerical solution of the transport equation},
  year        = {1976},
  month       = jul,
  doi         = {10.2172/104974},
  school      = {Office of Scientific and Technical Information (OSTI)},
}

@Article{levermore1996,
  author    = {Levermore, C. David},
  journal   = {J. Stat. Phys.},
  title     = {Moment closure hierarchies for kinetic theories},
  year      = {1996},
  issn      = {1572-9613},
  month     = jun,
  number    = {5–6},
  pages     = {1021--1065},
  volume    = {83},
  doi       = {10.1007/bf02179552},
  publisher = {Springer},
}

@Article{minerbo1978,
  author    = {Minerbo, Gerald N.},
  journal   = {J. Quant. Spectrosc. Radiat. Transfer},
  title     = {Maximum entropy {E}ddington factors},
  year      = {1978},
  issn      = {0022-4073},
  month     = dec,
  number    = {6},
  pages     = {541--545},
  volume    = {20},
  doi       = {10.1016/0022-4073(78)90024-9},
  publisher = {Elsevier},
}

@PhdThesis{monreal2013,
  author      = {Monreal, Philipp and Frank, Martin},
  school      = {{RWTH} Aachen University},
  title       = {{Moment Realizability and Kershaw Closures in Radiative Transfer}},
  year        = {2013},
  institution = {Lehr-und Forschungsgebiet Simulation in der Kerntechnik},
  url         = {http://publications.rwth-aachen.de/record/210538},
}

@Article{moujaes2025,
  author    = {Moujaes, Paul and Kuzmin, Dmitri},
  journal   = {J. Comput. Phys.},
  title     = {Monolithic convex limiting and implicit pseudo-time stepping for calculating steady-state solutions of the {E}uler equations},
  year      = {2025},
  pages     = {113687},
  volume    = {523},
  doi       = {10.1016/j.jcp.2024.113687},
  publisher = {Elsevier},
}

@Article{olbrant2012,
  author    = {Olbrant, Edgar and Hauck, Cory D. and Frank, Martin},
  journal   = {J. Comput. Phys.},
  title     = {A realizability-preserving discontinuous {G}alerkin method for the {M}1 model of radiative transfer},
  year      = {2012},
  number    = {17},
  pages     = {5612--5639},
  volume    = {231},
  doi       = {10.1016/j.jcp.2012.03.002},
  publisher = {Elsevier},
}

@Article{pichard2016,
  author    = {Pichard, Teddy and Aregba-Driollet, Denise and Brull, St{\'e}phane and Dubroca, Bruno and Frank, Martin},
  journal   = {Commun. Comput. Phys.},
  title     = {Relaxation schemes for the {$M_1$} model with space-dependent flux: {A}pplication to radiotherapy dose calculation},
  year      = {2016},
  number    = {1},
  pages     = {168--191},
  volume    = {19},
  doi       = {10.4208/cicp.121114.210415a},
  publisher = {Cambridge University Press},
}

@Article{pichard2017,
  author    = {Pichard, Teddy and Alldredge, Graham W. and Brull, St{\'e}phane and Dubroca, Bruno and Frank, Martin},
  journal   = {J. Sci. Comput.},
  title     = {An approximation of the {$M_2$} closure: {A}pplication to radiotherapy dose simulation},
  year      = {2017},
  number    = {1},
  pages     = {71--108},
  volume    = {71},
  doi       = {10.1007/s10915-016-0292-8},
  publisher = {Springer},
}

@Article{frank2007,
  author    = {Frank, Martin and Hensel, Hartmut and Klar, Axel},
  journal   = {SIAM J. Appl. Math.},
  title     = {A fast and accurate moment method for the {F}okker–{P}lanck equation and applications to electron radiotherapy},
  year      = {2007},
  issn      = {1095-712X},
  month     = jan,
  number    = {2},
  pages     = {582--603},
  volume    = {67},
  doi       = {10.1137/06065547x},
  publisher = {Society for Industrial & Applied Mathematics (SIAM)},
}

@Article{Levermore1984,
  author    = {Levermore, C. David},
  journal   = {J. Quant. Spectrosc. Radiat. Transfer},
  title     = {Relating {E}ddington factors to flux limiters},
  year      = {1984},
  issn      = {0022-4073},
  month     = feb,
  number    = {2},
  pages     = {149--160},
  volume    = {31},
  doi       = {10.1016/0022-4073(84)90112-2},
  publisher = {Elsevier BV},
}

@Article{bedford2019,
  author    = {Bedford, James L.},
  journal   = {Phys. Med. Biol.},
  title     = {Calculation of absorbed dose in radiotherapy by solution of the linear {B}oltzmann transport equations},
  year      = {2019},
  number    = {2},
  pages     = {02TR01},
  volume    = {64},
  doi       = {10.1088/1361-6560/aaf0e2},
  publisher = {IOP Publishing},
}

@Article{vassiliev2010,
  author    = {Vassiliev, Oleg N. and Wareing, Todd A. and McGhee, John and Failla, Gregory and Salehpour, Mohammad R. and Mourtada, Firas},
  journal   = {Phys. Med. Biol.},
  title     = {Validation of a new grid-based {B}oltzmann equation solver for dose calculation in radiotherapy with photon beams},
  year      = {2010},
  number    = {3},
  pages     = {581},
  volume    = {55},
  publisher = {IOP Publishing},
}

@Article{gifford2006,
  author    = {Gifford, Kent A. and Horton, John L. and Wareing, Todd A. and Failla, Gregory and Mourtada, Firas},
  journal   = {Phys. Med. Biol.},
  title     = {Comparison of a finite-element multigroup discrete-ordinates code with {M}onte {C}arlo for radiotherapy calculations},
  year      = {2006},
  number    = {9},
  pages     = {2253},
  volume    = {51},
  doi       = {10.1088/0031-9155/51/9/010},
  publisher = {IOP Publishing},
}

@Article{verbeek2021,
  author    = {Verbeek, Nico and Wulff, J{\"o}rg and Janson, Martin and B{\"a}umer, Christian and Zahid, Sameera and Timmermann, Beate and Brualla, Lorenzo},
  journal   = {Med. Phys.},
  title     = {Experiments and {M}onte {C}arlo simulations on multiple {C}oulomb scattering of protons},
  year      = {2021},
  issn      = {2473-4209},
  month     = may,
  number    = {6},
  pages     = {3186--3199},
  volume    = {48},
  doi       = {10.1002/mp.14860},
  publisher = {Wiley Online Library},
}

@Article{Saini2018,
  author   = {Jatinder Saini and Erik Traneus and Dominic Maes and Rajesh Regmi and Stephen R. Bowen and Charles Bloch and Tony Wong},
  journal  = {Translational Lung Cancer Research},
  title    = {Advanced proton beam dosimetry part {I}: {R}eview and performance evaluation of dose calculation algorithms},
  year     = {2018},
  issn     = {2226-4477},
  number   = {2},
  volume   = {7},
  url      = {https://tlcr.amegroups.com/article/view/20846},
}

@Article{Lin2021,
  author    = {Lin, Liyong and Taylor, Paige A. and Shen, Jiajian and Saini, Jatinder and Kang, Minglei and Simone II, Charles B. and Bradley, Jeffrey D. and Li, Zuofeng and Xiao, Ying},
  journal   = {International Journal of Particle Therapy},
  title     = {{NRG} oncology survey of {Monte Carlo} dose calculation use in {US} proton therapy centers},
  year      = {2021},
  number    = {2},
  pages     = {73--81},
  volume    = {8},
  doi       = {10.14338/ijpt-d-21-00004},
  publisher = {Elsevier},
}

@Article{JANSON2024,
  author  = {Martin Janson and Lars Glimelius and Albin Fredriksson and Erik Traneus and Erik Engwall},
  journal = {Med. Dosim.},
  title   = {Treatment planning of scanned proton beams in {RayStation}},
  year    = {2024},
  issn    = {0958-3947},
  note    = {Treatment Planning in Proton Therapy},
  number  = {1},
  pages   = {2-12},
  volume  = {49},
  doi     = {https://doi.org/10.1016/j.meddos.2023.10.009},
  url     = {https://www.sciencedirect.com/science/article/pii/S095839472300105X},
}

@Article{duclous2010,
  author    = {Duclous, Roland and Dubroca, Bruno and Frank, Martin},
  journal   = {Phys. Med. Biol.},
  title     = {A deterministic partial differential equation model for dose calculation in electron radiotherapy},
  year      = {2010},
  number    = {13},
  pages     = {3843},
  volume    = {55},
  publisher = {IOP Publishing},
}

@Article{Gottschalk2009,
  author    = {Gottschalk, Bernard},
  journal   = {Medical Physics},
  title     = {On the scattering power of radiotherapy protons},
  year      = {2009},
  issn      = {2473-4209},
  month     = dec,
  number    = {1},
  pages     = {352--367},
  volume    = {37},
  doi       = {10.1118/1.3264177},
  publisher = {Wiley},
}

@Misc{Stammer2025,
  author    = {Stammer, Pia and Wahl, Niklas and Kusch, Jonas and Lathouwers, Danny},
  title     = {A high-order deterministic dynamical low-rank method for proton transport in heterogeneous media},
  year      = {2025},
  copyright = {Creative Commons Attribution 4.0 International},
  doi       = {10.48550/ARXIV.2508.04484},
  publisher = {arXiv},
}

@Book{IAEA.2024,
  author    = {null IAEA.},
  publisher = {International Atomic Energy Agency},
  title     = {{Neutron Monitoring for Radiation Protection}},
  year      = {2024},
  address   = {Vienna},
  edition   = {1st ed.},
  isbn      = {9789201512222},
  note      = {Description based on publisher supplied metadata and other sources.},
  number    = {v.115},
  series    = {Safety Reports Series No Series},
  pagetotal = {1254},
  ppn_gvk   = {1908948183},
}

@Article{ashby2025,
  author    = {Ashby, Ben S. and Chronholm, Veronika and Hajnal, Daniel K. and Lukyanov, Alex and MacKenzie, Katherine and Pim, Aaron and Pryer, Tristan},
  journal   = {J. Math. Biol.},
  title     = {Efficient proton transport modelling for proton beam therapy and biological quantification},
  year      = {2025},
  issn      = {1432-1416},
  month     = apr,
  number    = {5},
  volume    = {90},
  doi       = {10.1007/s00285-025-02212-1},
  publisher = {Springer Science and Business Media LLC},
}

@Article{ulmer2007,
  author    = {Ulmer, W.},
  journal   = {Radiat. Phys. Chem.},
  title     = {Theoretical aspects of energy–range relations, stopping power and energy straggling of protons},
  year      = {2007},
  issn      = {0969-806X},
  month     = jul,
  number    = {7},
  pages     = {1089--1107},
  volume    = {76},
  doi       = {10.1016/j.radphyschem.2007.02.083},
  publisher = {Elsevier BV},
}

@Article{cox2024,
  author    = {Cox, Alexander M. G. and Hattam, Laura and Kyprianou, Andreas E. and Pryer, Tristan},
  journal   = {Proceedings of the Royal Society A: Mathematical, Physical and Engineering Sciences},
  title     = {A {B}yesian inverse approach to proton therapy dose delivery verification},
  year      = {2024},
  issn      = {1471-2946},
  month     = nov,
  number    = {2301},
  volume    = {480},
  doi       = {10.1098/rspa.2023.0836},
  publisher = {The Royal Society},
}

@Article{bortfeld1997,
  author    = {Bortfeld, Thomas},
  journal   = {Med. Phys.},
  title     = {An analytical approximation of the {B}ragg curve for therapeutic proton beams},
  year      = {1997},
  issn      = {2473-4209},
  month     = dec,
  number    = {12},
  pages     = {2024--2033},
  volume    = {24},
  doi       = {10.1118/1.598116},
  publisher = {Wiley},
}

@Article{newhauser2015,
  author    = {Newhauser, Wayne D. and Zhang, Rui},
  journal   = {Phys. Med. Biol.},
  title     = {The physics of proton therapy},
  year      = {2015},
  issn      = {1361-6560},
  month     = mar,
  number    = {8},
  pages     = {R155--R209},
  volume    = {60},
  doi       = {10.1088/0031-9155/60/8/r155},
  publisher = {IOP Publishing},
}

@Misc{berger1998,
  author    = {Berger, Martin J. and Coursey, J. S. and Zucker, M. A. and Chang, J.},
  title     = {Stopping-powers and range tables for electrons, protons, and helium Ions, {NIST Standard Reference Database 124}},
  year      = {1993},
  copyright = {License Information for NIST data},
  doi       = {10.18434/T4NC7P},
  language  = {en},
  publisher = {National Institute of Standards and Technology},
}

@Article{moujaes2026,
  author    = {Moujaes, Paul and Kuzmin, Dmitri and Bäumer, Christian},
  journal   = {Math. Comput. Simulat.},
  title     = {Realizability-preserving monolithic convex limiting in continuous {G}alerkin discretizations of the {M}1 model of radiative transfer},
  year      = {2026},
  issn      = {0378-4754},
  month     = aug,
  pages     = {570--590},
  volume    = {246},
  doi       = {10.1016/j.matcom.2026.02.011},
  publisher = {Elsevier BV},
}

@Article{morel1981,
  author    = {Morel, J. E.},
  journal   = {Nucl. Sci. Eng.},
  title     = {{F}okker-{P}lanck calculations using standard discrete ordinates transport codes},
  year      = {1981},
  issn      = {1943-748X},
  month     = dec,
  number    = {4},
  pages     = {340--356},
  volume    = {79},
  doi       = {10.13182/nse79-340},
  publisher = {Informa UK Limited},
}

@Article{hensel2006,
  author    = {Hensel, Hartmut and Iza-Teran, Rodrigo and Siedow, Norbert},
  journal   = {Phys. Med. Biol.},
  title     = {Deterministic model for dose calculation in photon radiotherapy},
  year      = {2006},
  issn      = {1361-6560},
  month     = jan,
  number    = {3},
  pages     = {675--693},
  volume    = {51},
  doi       = {10.1088/0031-9155/51/3/013},
  publisher = {IOP Publishing},
}

@MastersThesis{ulikema2012,
  author = {Uilkema, S. B.},
  school = {TU Delft},
  title  = {{Proton Therapy Planning using the SN Method with theFokker-Planck Approximation}},
  year   = {2012},
}

@Article{berthon2011,
  author    = {Berthon, Christophe and Frank, Martin and Sarazin, Céline and Turpault, Rodolphe},
  journal   = {Commun. Comput. Phys.},
  title     = {Numerical methods for balance Laws with space dependent flux: {A}pplication to radiotherapy dose calculation},
  year      = {2011},
  issn      = {1991-7120},
  month     = nov,
  number    = {5},
  pages     = {1184--1210},
  volume    = {10},
  doi       = {10.4208/cicp.020810.171210a},
  publisher = {Global Science Press},
}

@Article{larsen1997,
  author    = {Larsen, Edward W. and Miften, Moyed M. and Fraass, Benedick A. and Bruinvis, Ia{\i}n A. D.},
  journal   = {Med. Phys.},
  title     = {{Electron dose calculations using the Method of Moments}},
  year      = {1997},
  issn      = {2473-4209},
  month     = jan,
  number    = {1},
  pages     = {111--125},
  volume    = {24},
  doi       = {10.1118/1.597920},
  publisher = {Wiley},
}

\end{document}